\documentclass{article}

\usepackage{amsmath}
\usepackage{amsfonts}
\usepackage{amssymb}
\usepackage{color}  
\usepackage{amsthm}
\usepackage{graphicx}
\usepackage{caption}
\usepackage{comment}
\usepackage{subcaption}
\usepackage{algpseudocode}
\usepackage{algorithm}
\usepackage[utf8]{inputenc}
\usepackage{tabularx}
\usepackage{multirow}
\usepackage{hyperref}
\usepackage[english]{babel}
\usepackage{setspace}

\usepackage{natbib}
\bibpunct[, ]{(}{)}{,}{a}{}{,}%

\title{Local Area Routes and Valid Inequalities for Efficient Vehicle Routing%:  A principled synthesis of Graph Generation, Local Area Routes and Subset Row Inequalities
}

\author{Udayan Mandal\textsuperscript{ \rm 1,\rm 2}, Amelia Regan\textsuperscript{\rm 2, \rm 3}, Julian Yarkony \textsuperscript{\rm 4} \\% All authors must be in the same font size and format. Use \Large and \textbf to achieve this result when breaking a line
\textsuperscript{\rm 1}Stanford University, Palo Alto, CA\\
\textsuperscript{\rm 2}University of California, Irvine, CA\\
\textsuperscript{\rm 3}University of Washington, Seattle, WA\\
\textsuperscript{\rm 4} Laminaar Optimization Research Group, La Jolla, CA
}\date{September, 26, 2022 with update on January, 30, 2023}
\begin{document}
\maketitle

\begin{abstract}
    
In this research we introduce Local Area (LA) routes for improving the efficiency and tightness of column generation (CG) methods for solving vehicle routing problems (VRP).
LA-routes rely on pre-computing the lowest cost elementary sub-route (called an LA-arc) for each tuple consisting of the following: \textbf{(1)} a (first) customer where the LA-arc begins, \textbf{(2)} a distant customer (from the first) where the LA-arc ends, and \textbf{(3)} a set of intermediate customers near the first customer. LA-routes are constructed by concatenating LA-arcs where the final customer in a given LA-arc is the first customer in the subsequent LA-arc. A Decremental State Space Relaxation method is applied over LA-routes to construct the lowest reduced cost elementary route during the pricing step of CG. LA-route based solvers can be used to efficiently tighten the standard set cover VRP using a variant of subset row inequalities, which do not alter the structure of pricing. 

We incorporate LA-arcs into a novel CG stabilization scheme.  Specifically each column generated during pricing is mapped to an ordered list of customers consistent with that column.  An LA-arc is consistent with an ordering if the first/last customer in the arc come before/after all other customers in the LA-arc in the associated ordering respectively.  Each such ordering is then mapped to a multi-graph where nodes correspond to (customer/demand) and edges correspond to LA-arcs consistent with that ordering. Hence any path from source to sink on the multi-graph is a feasible elementary route.  The ordering for a column places customers spatially nearby in nearby positions on the ordering so that routes can be generated so as to permit spatially nearby customers to be visited without traveling far away first. We solve the restricted master problem over these graphs, which has special structure allowing for fast solution.
\end{abstract}

%\KEYWORDS{deterministic inventory theory; infinite linear programming duality; 
%  existence of optimal policies; semi-Markov decision process; cyclic schedule}

%\KEYWORDS{column generation; CVRP; cut and price; route relaxation}

\maketitle
%%%%%%%%%%%%%%%%%%%%%%%%%%%%%%%%%%%%%%%%%%%%%%%%%%%%%%%%%%%%%%%%%%%%%%

% Samples of sectioning (and labeling) in TRSC
% NOTE: (1) \section and \subsection do NOT end with a period
%       (2) \subsubsection and lower need end punctuation
%       (3) capitalization is as shown (title style).
%
%\section{Introduction.}\label{intro} %%1.
%\subsection{Duality and the Classical EOQ Problem.}\label{class-EOQ} %% 1.1.
%\subsection{Outline.}\label{outline1} %% 1.2.
%\subsubsection{Cyclic Schedules for the General Deterministic SMDP.}
%  \label{cyclic-schedules} %% 1.2.1
%\section{Problem Description.}\label{problemdescription} %% 2.

% Text of your paper here

%\include{appendix.tex}
\section{Introduction}
In this paper we consider an approach for improving the efficiency of column generation (CG)\citep{barnprice,cuttingstock} methods for solving vehicle routing problems (VRP). Such problems have had rich applications in industrial, military and humanitarian logistics for decades, and emerging problems including rapid grocery delivery, the routing of robots in warehouses, and routing of autonomous vehicles on street networks have made these problems especially interesting.  %In this paper we consider develop a new algorithmic tool for tackling problems in vehicle routing called local area routes
We introduce Local Area (LA) route relaxations, an alternative/complement to the commonly used ng (\textbf{n}ei\textbf{g}hborhood)-route relaxations \citep{baldacci2011new}, and Decremental State Space Relaxations (DSSR) \citep{righini2008new} inside of CG formulations. LA-routes are a subset of ng-routes and a super-set of elementary routes. Normally, the pricing stage of CG must produce elementary routes, which are routes without repeated customers, using processes which can be computationally expensive. Non-elementary routes visit at least one customer more than once, creating a cycle. LA-routes relax the constraint of being an elementary route in such a manner as to permit efficient pricing.  Local Area (LA) routes, can serve as a component in CG solutions \citep{barnprice,cuttingstock} to VRP \citep{Desrochers1992,feillet2010tutorial}. While our approach can be adapted/applied to general VRP, we describe and conduct experiments using the Capacitated Vehicle Routing Problem (CVRP), which we describe using the following terms: \vspace{2mm} \\\textbf{(a)} A depot located in space \\\textbf{(b)} A set of customers with integer demands located in space  \\\textbf{(c)} A set of homogeneous vehicles with integer capacity.\vspace{2mm} \\  Vehicles are assigned to routes, where each route satisfies the following: \\\textbf{(1)} Each route starts and ends at the depot\\ \textbf{(2)} The total demand of the customers on a route does not exceed the capacity of a vehicle\\ \textbf{(3)} The cost of the route is the total distance traveled\\ \textbf{(4)} The vehicle leaves a customer the same number of times as it arrives at a customer (no teleportation of vehicles)\\ \textbf{(5)} No customer is serviced more than once on a route \\The CVRP problem selects a set of routes with the goal of minimizing the total distance traveled while ensuring that each customer is serviced at least once.  As an aside for purposes of later notation, if a route satisfies properties \textbf{(1,2,3,4)} but not necessarily \textbf{(5)} then it is referred to as non-elementary.   % but otherwise feasible.

CVRP and other VRP can be solved with compact linear programming (LP) relaxations or expanded LP relaxations \citep{Desrochers1992} (where the expanded LP relaxation is often referred to as the set cover formulation). Compact LP relaxations (CLP) have a variable for each tuple consisting of a vehicle and pair of destinations, each either at a depot or customer, denoting the first and second destination. This variable is used to indicate if the given vehicle travels from the first destination to the second destination. Expanded LP relaxations (ELP) have one variable for each possible route, and each such variable is used to indicate that the given route is selected. While both the CLP and ELP have the same optimal integer linear programming (ILP) objective function, these formulations can have different optimal LP objective values. The ELP is often tighter \citep{geoffrion1974lagrangean}, and is never looser than the CLP, leading it to be preferred in practice \citep{costa2019}. We should note that while no optimal integer solution to the ELP would use a route covering a customer more than once, this is not necessarily true for optimal fractional solutions to the ELP where fractions of a route covering a customer more than once can be used. Hence enforcing property \textbf{(5)} of routes is crucial for the ELP.  %preferred in practice to the compact relaxation since it tends to have a significantly larger optimal  LP value than the compact relaxation , meaning that the expanded relaxation is meaningfully tighter than the compact formulation.
The number of possible routes can grow rapidly in the number of customers, so it is often infeasible to enumerate all routes. Hence the ELP is solved using CG, which imitates the revised simplex, where the pricing operation is a resource constrained shortest path (RCSP) problem \citep{baldacci2011new,Desrochers1992,irnich2005shortest}. The RCSP is NP-hard but is efficiently solvable at the scale of many practical problems \citep{desrosiers2005primer}.

We seek to ease solving the pricing problem by producing a class of routes known as Local Area (LA) routes that are easy to price over, and can be used to efficiently generate elementary routes when used inside the DSSR \citep{righini2008new,righini2009decremental}.  
An LA-route can be a non-elementary route, and must not contain cycles localized in space. Note that a cycle is a section of a route consisting of the same customer at the start and end of the section. Localized cycles in space are cycles consisting of customers that are all spatially close to one another. LA-routes are a tighter version of the celebrated ng-routes \citep{baldacci2011new}. LA-routes also naturally admit a stabilized formulation where stabilization alters the path of dual solutions so as to converge to the optimum faster \citep{du1999stabilized,marsten1975boxstep,ben2006dual,rousseau2007interior, oukil2007stabilized, haghani2020smooth}. We present our computational results in Section \ref{sec_exper}. To captivate our readers' attention, we discuss a few findings here. Our codes are experimental and far from production level, so the important details are the relative solution times and relative number of iterations required to find optimal solutions. Using a commonly available set of benchmark instances (Augerat A) \citep{augerat1995computational} we %found that permitting LA-routes of sizes 6-10 
typically speed up pricing times by a factor of about 60, and reduce the number of iterations required by as much as a factor of 200 (typically also a factor of about 60). Overall speedups are typically more than 10 times that of the baseline CG, though the additional time is consumed by overhead, and efficient coding in an industrial version of these codes would be much faster, thus making the computation time look more like the pricing time reduction.% By efficiently coding and parallelizing the generation of LA-arcs  we expect to be able to make this step negligible in time making the overall time speed up match the pricing time speedup.  %and  we can achieve improve this overall time considerably, perhaps to match the reductions in pricing time and iterations required.
When we refer to a baseline CG here we mean a basic column generation with no additional stabilization. 

We organize this document as follows. In Section \ref{Sec_math_rev} we provide a mathematical review of CG for CVRP. In Section \ref{litRev} we review the literature most related to this work. In Section \ref{sec_LA} we describe LA-routes, and consider their use in the CG pricing problem solving for the lowest reduced cost elementary route. In Section \ref{Sec_SRI_LA} we introduce LA-subset row inequalities  (LA-SRI)  (which are an adaptation of SRI \citep{jepsen2008subset}) to tighten the ELP, and show how the structure of the pricing problem is preserved when considering LA-SRI.  %IN Section \ref{Sec_GG} we show how to apply GG in such a manner as to use the SRI.
In Section \ref{sec_new_gg_pgm} we describe an efficient stabilized solver for the ELP that is able to consider LA-SRI. In Section \ref{sec_exper} we provide an experimental validation of our approach. In Section \ref{sec_conc} we conclude and discuss extensions and we provide additional mathematical exposition along with experimental results in the Appendix. Before we move on we present relevant notation in Table \ref{basicTab}.
\newpage
\begingroup
\setlength{\tabcolsep}{10pt} % Default value: 6pt
\renewcommand{\arraystretch}{1.5} % Default value: 1
\begin {table}[ht]
\begin{center}
\begin{tabular}{|p{1.2cm}|p{3cm}||p{10cm}|}
 \hline
Name & Space & Meaning\\
 \hline
   $N$ & Set & Set of customers\\
   $N^+$ & Set & Set of customers plus the depot which is counted as $-1$ for starting and $-2$ for ending even though they are the same place\\
$d_0$ & Scalar & Amount of capacity in a vehicle \\
 $d_u$ & Scalar & Amount demand at customer $u$ \\
$c_{uv}$ & $u \in N^+,v \in N^+$&  Distance from $u$ to $v$.\\
$\pi_u$ & $u \in N$ & Dual variable associated with customer $u$\\
$\pi_0$ & Scalar & Dual variable associated with enforcing an upper bound on the number of vehicles used  \\
$\Omega$& Set & Set of elementary routes.  \\
$\Omega^0$ & Set & Set of non-elementary routes.  Routes in $\Omega^0$ satisfy properties (1,2,3,4) from the intro but not necessarily property $5$.  Thus  $\Omega \subseteq \Omega^0$.\\% that obey all properties except elementarity
$a_{ul}$ & $u \in N,l \in \Omega^0$ & $a_{ul}$ is the number of times route $l$ covers customer $u$.  Note that $a_{ul}$ is binary for $l \in \Omega$ but not necessarily for $l\in \Omega^0$.\\
$a_{uvdl}$ & $u \in N^+-(-2),v \in N^+-(u\cup -1),d \in D,l \in \Omega^0$ & $a_{uvdl}=1$ if the route $l$ leaves $u$ with $d$ units of capacity remaining and travels immediately to $v$. Here $u,v,d$ must also lie in $F$ which is the set of feasible possibilities of $u,v,d$. \\
 \hline
\end{tabular}
\end{center}
\vspace{0.25cm}
\caption {\bf{CVRP Route Notation}} 
\label{basicTab}
\end{table}
\endgroup
 
\section{Mathematical Review}
\label{Sec_math_rev}
We now describe CVRP formally using the following notation. We use $N$ to denote the set of customers. We use $N^+$ to denote $N$ augmented with the starting/ending depot denoted $-1,-2$ respectively, which are co-located. We use $K$ to denote the maximum number of vehicles that can be used, where each vehicle has capacity $d_0$. % to denote the capacity of a vehicle.  
We use $d_{u}$ to denote the demand of $u \in N^+$ where $d_{-1}=d_{-2}=0$. We use $\Omega$ to denote the set of feasible routes, which we index by $l$. The set of non-elementary routes is denoted $\Omega^0$.  We set $a_{ul}=1$ if route $l$ services customer $u$. We set $a_{uvdl}=1$ if route $l$ leaves $u \in N^+-(-2)$ with $d$ units of capacity remaining then travels immediately to $v \in N^+-(-1\cup u)$.  We use $F$ to denote the set of feasible values $uvd$ for $a_{uvdl}$ which is defined for tuples $u\in N^+-(-2), v\in (N^+-(-1 \cup u)),(d_0-d_u)\geq d \geq d_v$.  
We use $c_{uv}$ to denote the distance between any pair $u\in N^+, v \in N^+$. Using $c_{uv}$, we write the cost of a route $l$ which we denote as $c_l$ as the total distance traveled with this equation: $c_l=\sum_{uvd\in F}c_{uv}a_{uvdl} \quad \forall l \in \Omega$.  
We now consider optimization for CVRP below as a weighted set cover problem which is the ELP. %We define $\Omega$ to be the set of all possible routes. 
The standard CVRP ELP is written using decision variables $\theta_l$ where $\theta_l=1$ if route $l$ (referred to as a column in CG) is selected in our solution, otherwise $\theta_l=0$. 
We write the ELP over $\Omega$ as $ \Psi(\Omega)$ below %(which is also referred to as the master problem (MP)),
with dual variables $\pi$ written in [].% as $\pi$. 
\begin{subequations}
\label{primal_master}
\begin{align}
    \Psi(\Omega)=\min_{\theta \geq 0}\sum_{l \in \Omega}c_l \theta_l \label{CVRP_obj}\\
    \sum_{l \in \Omega}a_{ul}\theta_l\geq 1 \quad \forall u \in N \quad [\pi_u] \label{CVRP_cover}\\
    \sum_{l \in \Omega}\theta_l\leq K \quad [-\pi_0] \label{CVRP_packing}
\end{align}
\end{subequations}
In \eqref{CVRP_obj} we minimize the total cost of the routes used.  In \eqref{CVRP_cover} we ensure that each customer is serviced at least once (though an optimal solution will service each customer exactly once). In \eqref{CVRP_packing} we ensure that no more than $K$ vehicles are used in the solution.  Since the size of $\Omega$ can grow exponentially in the number of customers we cannot trivially solve \eqref{primal_master} which is called the master problem (MP). Instead CG \citep{barnprice,Desrochers1992} is employed to solve \eqref{primal_master}. CG constructs a sufficient subset of $\Omega$ denoted $\Omega_R$ s.t. solving \eqref{primal_master} using $\Omega_R$ (denoted $\Psi(\Omega_R)$) provides an optimal solution to \eqref{primal_master} using $\Omega$. 
To construct $\Omega_R$, we iterate between \textbf{(1)} solving $\Psi(\Omega_R)$, which is referred to as the restricted master problem (RMP)
and \textbf{(2)} identifying at least one $l \in \Omega$ with negative reduced cost, which are then added to $\Omega_R$.
In many CG applications only the lowest reduced cost column in $ \Omega$ is generated. We write the selection of this route as optimization below using $\bar{c}_l$ to denote the reduced cost of route $l$.% \in \Omega$. 
\begin{subequations}
\label{pricing}
\begin{align}
    \min_{l \in \Omega} \bar{c}_l \\
    \bar{c}_l=c_l+\pi_0-\sum_{u \in N}a_{ul}\pi_u  \quad \forall l \in \Omega \label{redCostForm}
\end{align}
\end{subequations}
% %, where r).  %One or more negative reduced cost columns are generated during pricing.
The operation in \eqref{pricing} is referred to as pricing and solved as a resource constrained shortest path problem (RCSP)\citep{irnich2005shortest}, not by enumerating all of $\Omega$, which would be computationally prohibitive.  CG terminates when pricing proves no column with negative reduced cost exists in $\Omega$.  This certifies that CG has produced the optimal solution to \eqref{primal_master}. %We terminate CG when no $l$ in $\Omega$ has negative reduced cost.  
In CG $\Omega_R$ is typically initialized with artificial variables that have prohibitively high cost but form a feasible solution.  % This can be done by creating $|N|$ variables each of which covers a customer $u$ with prohibitively high cost but without using a vehicle.% Pricing is attacked as an elementary resource constrained shortest path problem (RCSP) and may be treated with a variety of algorithms including but not limited to dynamic programming based labeling algorithms\citep{desrosiers2005primer}.  
%\subsection{Subset Row Inequalities for CVRP}
%\label{Sec_SRI_basic}
The ELP can be tightened using valid inequalities such as subset row inequalities (SRI) \citep{jepsen2008subset}. The simplest example of SRI for CVRP is written as follows. For any given set of three unique customers denoted $N_{\delta} \subseteq N$ the number of routes including $2$ or more members of $N_{\delta}$ can not exceed one. 

We motivate SRI by the following example. Consider that we have a depot in San Diego and three customers nearly co-located in and lying in the distant city New York City (NYC) each with demand $1$ and vehicles of capacity $2$. The optimal ILP solution requires $2$ cross country routes. However the optimal fractional solution uses $\frac{3}{2}$ cross country routes with $\frac{1}{2}$ of each of the $3$ routes servicing exactly $2$ customers. Specifically each of the following routes are used with $\theta$ value equal to $\frac{1}{2}$:  (-1, NYC1, NYC2, -2),(-1, NYC1, NYC3, -2),(-1, NYC2, NYC3, -2); where the sequence describes the locations visited in order, and NYC1, NYC2, NYC3 indicate the three separate customers in NYC.   When enforcing that no more than one route servicing $2$ or more of the $3$ customers is selected then the LP becomes tight. Observe that the solution previously shown uses $\frac{3}{2}$ routes servicing $2$ of the $3$ customers. Other SRI can express that the number of routes containing $3$ or more customers from a group of $5$ customers can not exceed $1$. Similarly for any group of $5$ customers the number of routes servicing $2,3$ customers plus $2$ times the number of routes servicing $4,5$ customers can not exceed $2$. The general version of SRI states that for any set of customers $N_{\delta}$ and positive integer $m_{\delta}$, the following holds: 
\begin{align}
\label{my_div_2}
\sum_{l\in \Omega}\theta_l \left\lfloor\frac{\sum_{u \in N_{\delta}}a_{ul}}{m_{\delta}}\right\rfloor\leq \left\lfloor \frac{|N_{\delta}|}{m_{\delta}}  \right\rfloor \; \forall N_{\delta}\subseteq N,  m_{\delta} \in \mathbb{Z}_{+}
\end{align}
When \eqref{my_div_2} is added to the MP we solve the MP in a cutting plane manner where given any subset of valid inequalities we solve optimization using CG. The nascent set of such constraints is referred to as $\Delta_R$ and the associated RMP as $\Psi(\Omega_R,\Delta_R)$. The addition of such SRI to the RMP alters the structure of pricing, which can make pricing challenging as each SRI introduces a new resource that must be kept track of.  

\section{Literature Review}
\label{litRev}

\subsection{General Dual Stabilization}

The number of iterations of column generation (CG) required to optimally solve the master problem (MP) can be dramatically decreased by intelligently altering the sequence of dual solutions generated \citep{du1999stabilized, rousseau2007interior, Pessoa2018Automation} over the course of CG. Such approaches, called dual stabilization, can be written as seeking to maximize the Lagrangian bound at each iteration of CG \citep{geoffrion1974lagrangean}.  The Lagrangian bound is a lower bound on the optimal solution objective  to the MP that can be easily generated at each iteration of CG. In CVRP problems the Lagrangian bound is the LP value of the restricted master problem (RMP) plus the reduced cost of the lowest reduced cost column times the number of customers. Observe that when no negative reduced cost columns exist, the Lagrangian bound is simply the LP value of the RMP. The Lagrangian bound is a concave function of the dual variable vector. The current columns in the RMP provide for a good approximation of the Lagrangian bound nearby dual solutions generated thus far but not regarding distant dual solutions.  This motivates the idea of attacking the maximization of the Lagrangian bound in a manner akin to gradient ascent. Specifically we trade off maximizing the objective of the RMP, and having the produced dual solution be close to the dual solution with the greatest Lagrangian bound identified thus far (called the incumbent solution).  

A simple but effective version of this idea is the the box-step method of \citep{marsten1975boxstep}, which maximizes the Lagrangian bound at each iteration of CG s.t. the dual solution does not leave a bounding box around the incumbent solution.  Given the new solution, the lowest reduced cost column is generated and if the associated Lagrangian bound is greater than that of the incumbent then the incumbent is updated. The simple approach of \citep{Pessoa2018Automation} takes the weighted combination of the incumbent solution and the solution to the RMP and performs pricing on that weighted combination.  %dual stabilization can make CG behave like a gradient ascent method 
%Due to inherent instability of the dual variables in CG, many methods of stabilization have been proposed over the last two decades.
Du Merle et al formalized the idea of stabilized CG in their 1999 paper of that name \citep{du1999stabilized}. That paper proposed a 3-piecewise linear penalty function to stabilize CG. Ben Amor and Desrosiers later proposed a 5-piecewise linear penalty function for improved stabilization \citep{amor2006proximal}. Shortly after, Oukil et al used the same framework to attack highly degenerate instances of multiple-depot vehicle scheduling problems \citep{oukil2007stabilized}. Ben Amor et al later proposed a general framework for stabilized CG algorithms in which a stability center is chosen as an estimate of the optimal dual solution \citep{amor2009choice}. Gonzio et al proposed a primal-dual CG method in which the sub-optimal solutions of the RMP are obtained using an interior point solver that is proposed in an earlier paper by the first author \citep{gondzio1995hopdm}. They examine their solution method relative to standard CG and the analytic center cutting plane method proposed by Babonneau et al. \citep{babonneau2006solving, babonneau2007proximal}. They found that while standard CG is efficient for small problem instances, the primal-dual CG method performed better than standard CG on larger problems \citep{Gondzio2013New}.% Our experimental analysis had the same findings. 

%\subsubsection{Trust Region Based Stabilization}
%Trust region based methods exploit the understanding of CG operating a search algorithm over the dual space \citep{marsten1975boxstep}. Such methods seek to maximize the Lagrangian Relaxation of the master problem. Since columns generated at a given dual solution do not provide good information regarding the Lagrangian relaxation at distant points in the dual space CG can have a problem of bouncing between distant points in the dual space. This is circumvented by establishing a trust region computed at the point corresponding to the greatest Lagrangian relaxation thus far identified, and limiting search around it \citep{marsten1975boxstep}. This approach is extend in \citep{du1999stabilized} in which dual variables are penalized for leaving the trust region. 

\subsection{Dual Optimal Inequalities}

Dual optimal inequalities (DOI) \citep{ben2006dual} provide easily computed provable bounds on the space where the optimal dual solution to the MP lies. In this manner, the use of DOI reduces the size of the dual space that CG must search over and hence the number of iterations of CG. Dual constraints corresponding to DOI are typically defined over one or a small number of variables and hence do not significantly increase the solution time of the RMP (though exceptions exist \citep{haghani2020smooth}). DOI are problem instance or problem domain specific. One such example of DOI is in problems such as CVRP or cutting stock where the cost of any column is not increased by removing customers (or rolls in the cutting stock problem \citep{cuttingstock,lubbecke2005selected}) from the column.  Thus equality constraints enforcing that each customer is covered at least once (for the CVRP example) can be replaced by inequality constraints since the cost of a route is not increased by removing a customer from the route. In the dual representation, this replacement corresponds to enforcing that the dual variable corresponding to the constraint that a customer must be covered is non-negative.  

For the cutting stock problem, we can swap a roll of higher length for one of lower length without altering the feasibility of a column (alternatively known as a pattern).  Thus, it can be established that the dual variables associated with rolls, when ordered by non-decreasing roll size must be non-decreasing \citep{ben2006dual}.  In the primal representation, these bounds correspond to swap operations permitting a roll of a given length to be swapped for one of smaller length.  

In \citep{haghani2020smooth}, it is observed that the improvement in the objective corresponding to removing a customer from a column in CVRP (and also the Single Source Capacitated Facility Location Problem (SSCFLP)) can be bounded. Thus primal operations corresponding to removing customers from columns are provided.  In the dual representation, these operations enforce that the reduced cost of a column should not trivially become negative if customers are removed from it.  In the case of SSCFLP, for a given column, this property states that the dual contribution to the reduced cost for a given customer (included in the column) is treated as the maximum of the following two values: the dual variable for that customer, or the distance from the customer to the facility of that column.  

In \citep{haghani2020smooth} it is observed that the dual variables associated with constraints in problems embedded in a metric space should change smoothly over that space.  This is because the dual variable associated with a given customer roughly describes how much larger the objective of the LP is as a consequence of the given customer existing.  Thus, nearby customers should not normally have vastly different dual variables.  Specifically  \citep{haghani2020smooth} shows that in CVRP for any given pair of customers $u,v$ where $u$ has demand no less than $v$, the dual variable of $u$ plus two times the distance from $u$ to $v$ is no less than the dual variable of $v$.  In the primal form, any such pair corresponds to slack variables that provide for the swap operation from $u$ to $v$.  % customers.
\section{Local Area Routes}
\label{sec_LA}
In this section we describe our Local Area (LA) route relaxation based solver for CVRP pricing problems. We organize this section as follows. In Section \ref{sec_various_routes} we discuss various classes of routes relevant to ng-\citep{baldacci2011new} and LA- route relaxations. In Section \ref{sec_pricing_graph} we define the notation and basic information associated with LA-routes.  %In Section \eqref{sec_la_routes_pricing} we provide the definition for LA-routes.
In Section \ref{sec_decremental_la} we describe the use of LA-routes inside of Decremental State Space Relaxation (DSSR)\citep{righini2008new} to produce the the lowest reduced cost route elementary route. %In Section \ref{sec_various_routes} we provide preliminary information related to the various types of routes that can be considered.
In Section \ref{sec_pricing_la} we describe the computation of the lowest reduced cost LA-route using the Bellman-Ford Algorithm with efficient computation of edge weights for the associated graph specified in Section \ref{subsubsec_LA_review}.  In Section \ref{astar_dssr} where we replace the Bellman-Ford algorithm with A* \citep{dechter1985generalized} for  computational efficiency. In Section \ref{subsec_dom} we alter the graph over which pricing is conducted slightly, so as to permit the efficient use of dominance criteria commonly used in RCSP algorithms.%  Appendix \ref{appendix_imp_details}.% \ref{appendix_imp_details}. In Section \eqref{review_intel_cycle} we describe the reduction in the space of LA-routes for DSSR so as to improve the efficiency of pricing.  

\subsection{Various Classes of Routes}
\label{sec_various_routes}
%LA-route relaxations are best understood in terms of ng-route relaxations. Ng (\textbf{N}ei\textbf{G}hborhood)-routes are are permitted to have non-localized cycles in space; this means that at least one intermediate customer in the cycle (hence referred to as a breaker) must be spatially far away from starting customer in the cycle. LA-routes are described using a set of special indexes corresponding to customers on the route ordered from the start to the end of the route. LA-route relaxations further restrict the set of permitted cycles beyond that of ng-routes by additionally enforcing that the breaker must be a located at a special index where the set of special indexes is defined recursively as follows. %This set of special indexes in a route is defined recursively from the start of the route, with
%The first special index in the route is at index 1 meaning that it is associated with the first customer in the route. The $k$'th special index corresponds to the first customer after the $k-1$'th special index, that is not a neighbor of (therefore, spatially far from) the customer located at the $k-1$'th special index. %
%We demonstrate that LA-route relaxations can significantly improve the computational speed of pricing when compared to the standard Decremental State Space Relaxation. 

\begin{itemize}
%\item \textbf{Elementary Routes:  }

%We define $\Omega$ to be the set of all elementary routes, which are paths in $\Omega^0$ (where $\Omega^0$ is set of routes that need not be elementary but are otherwise feasible. ) that visit no customer in the path more than once. 
%For reference we define $\Omega$ formally below as a subset of $\Omega^0$ using $|p|$ to denote the number of customers in the path $p$.
%An elementary route is one which in addition to being a path from the source to the sink, visits each customer included in the route no more than once. We define the set of elementary routes as follows. We define $\Omega$ to be the set of all elementary routes. %u^p_k needS A DEF, AS DEOS vP

%\begin{align}
%   \Omega=\{ p \in \Omega^0; u^p_{k_1}\neq u^p_{k_2} \quad \forall k_1,k_2 \quad s.t. (0<k_1<k_2\leq |p|) 
%\end{align}
%
%$\Omega$ consists of any path in $\Omega^0$ that does not visit the same more than once.
%
\item \textbf{Q-routes:  }

Q-routes are routes in $\Omega^0$ which can not have cycles of length $1$. This means that a route visiting customer 1, followed by customer 2, and then visiting customer 1 again is forbidden. However, a route visiting customer 1, then visiting customer 2, then visiting customer 3, and then visiting customer 1 after is not forbidden.  Note that $u$ can never be visited immediately after $u$ in any route so cycles of size $0$ are forbidden in $\Omega^0$.   Q-routes were introduced by Christofides, Mingozzi and Toth in order to aid CG in solving vehicle routing problems \citep{christofides1981exact}.  We use $u^l_{k}$ to denote the $k'th$ customer visited in route $l$.

We use $\Omega^1$ to denote the set of Q-Routes, which is defined as follows using $|l|$ to denote the number of customers visited in $l$. 
\begin{align}
    \Omega^1=\{ l \in \Omega^0; u^l_{k}\neq u^l_{k+2} \quad \forall k \quad s.t. (0<k, k+2\leq |l|) 
\end{align}
%\item \textbf{KQ-routes: }
%These include no cycles of length $<=$ K and a generalization of $Q$ routes
  Q-routes can be generalized as KQ-routes where KQ-routes enforce that a customer can only be visited again after visiting at least $K$ intermediate customers. Observe that Q-routes correspond to KQ-routes with K=1. We use the term KQ to draw the link to Q-routes, but this terminology is not standard in the literature.  

We use $\Omega^K$ to denote the set of KQ-Routes, which is defined as follows. 
\begin{align}
   \Omega^K=\{ l \in \Omega^0; u^l_{k_1}\neq u^l_{k_2} \quad \forall k_1,k_2 \quad s.t.  (0<k_1<k_2\leq \min(|l|,k_1+K)) 
\end{align}
 
\item \textbf{Ng-routes:  }

Ng-routes (\textbf{n}ei\textbf{g}hborhood) are highly celebrated and used by many researchers \citep{baldacci2011new}. Ng-routes are a subset of $\Omega^0$ that does not include routes containing  spatially localized cycles but permits non-spatially localized cycles.  Ng-routes rely on each customer being associated with a set of customers, which are close in proximity to that customer (also known as neighbors of that customer). This set of neighbors of $u$ is denoted $M_u$, where $u$ represents the customer the set is associated with. Ng-routes ban spatially localized cycles by enforcing that a cycle can only exist starting and ending at $u$ if there is an intermediate customer $v$ for which $u\notin M_v$.

We now formally specify the set of ng-routes denoted $\Omega^{ng}$. Any $l \in \Omega^0$ lies in $\Omega^{ng}$ if the following holds for all $|l|\geq k_2>k_1\geq 1$ (with exposition below).%for which .  %is written as follows given any $k_2>k_1$.
\begin{align}
\label{ng_deg2a}
    (u^l_{k_1}=u^l_{k_2}) \rightarrow \exists k_3, \mbox{s.t. } k_1<k_3<k_2, u^l_{k_1} \notin M_{u^l_{k_3}}
\end{align}
The premise of \eqref{ng_deg2a} (left hand side of the $\rightarrow $ in \eqref{ng_deg2a} ) is true if the same customer is at indexes $k_1$ and $k_2$.  The right hand side of \eqref{ng_deg2a} states that there must exist a customer $k_3$ that lies between $k_1$ and $k_2$ and does not consider $u^l_{k_1}$ to be a ng-neighbor.

%A route is feasible if for a given customer $u$,  $u$ can only be visited again if a node which does not consider $u$ as a neighbor is visited first. This property bans spatially localized cycles.% and in practice solves the same LP relaxation as elementary. 
%To reiterate, a route is an ng-route if for any cycle in the route that cycle contains a non-neighbor node that does not cosnisdier  first node in the cycle.    

\item \textbf{LA-routes:}

LA-routes are a subset of ng-routes but further restrict cycles. Thus the CVRP set cover LP relaxation over LA-routes is no looser than, and in fact potentially tighter than that over ng-routes. LA-routes are defined using LA neighborhood set $N_u$ (for each customer $u$) where $N_u$ consists of spatially nearby customers to $u$. The LA neighborhood sets $N_u$ are computationally easier to consider than $M_u$ and hence can be larger than $M_u$. LA-routes are defined with a set of special indexes associated with each path $p$.  % in addition to the standard indexes.  The standard index of the $k$'th customer visited in a route (chronologically) is $k$.  %being $k$ means that this is the $k$'th node visited after leaving the depot.  %The first node visited after the dep%Here an index is associated with each node in the route in order of visit from start to finish treating repeated customers as unique.%  Thus the $j'$ th customer visited in a route is associated with index $j$.%  A route starting/ending at the depot then visiting customers 1,2,3,1 in order would have customers 
%
%The set of special indexes is defined using the neighborhood set $N_u$  
The set of special indexes in a route $l$ is defined recursively from the start of the route, with the first special index being equal to one.  Let $q^l_j$ be the index of the $j$'th special index (meaning $q^l_1=1$).  %Given that  $u^p_{q^p_j}=v$, %, and corresponds to the first customer visited in the route. 
The $j$'th special index corresponds to the first customer after $q^l_{j-1}$ that is not considered to be the an LA neighbor of $u^l_{q^l_{j-1}}$. %We use $q_{j}^p$ to denote the index of the $j$'th special given a path $p$ to describe the 
We define the set of special indexes by defining the $q^l_j$ terms recursively as follows.  %recursive definition of the special indices of a path $p$. % where the set of such indexes is denoted $Q^p$.  %.  'th customer in $p$.
\begin{align}
\label{qpdef}
    q_{j}^l\leftarrow \min_{\substack{k>q_{j-1}^l\\ u^l_k \notin N_v}} k \quad \forall j>1; \mbox{where $u^l_{q_{j-1}}=v$}\\
    q_1^l=1 \nonumber 
\end{align}
%We show that this can significantly speed up pricing, thereby allowing us to solve larger problems in real-time.
%To reiterate, a route is an ng-route if for any cycle in the route that cycle contains a non-neighbor node.
Any $l\in \Omega^0$ is an LA-route if for any cycle in that path starting/ending at customer $u$, there is an intermediate special index with associated customer $v$ for which $u \notin M_v$ (note the use of $M_v$ not $N_v$ here).  Note that the difference between an ng-route and an LA-route is understood by building on the statement ``Ng-routes ban spatially localized cycles by enforcing that a cycle can only exist starting and ending at $u$ if there is an intermediate customer $v$ for which $u\notin M_v$."  In contrast LA-routes ban spatially localized cycles by enforcing that a cycle can only exist starting and ending at $u$ if there is an intermediate customer $v$ {\textbf{at a special index}} for which $u\notin M_v$. We use $Q^l$ to refer to the set of special indexes in route $l\in \Omega^0$.  %Note that .   % Formally we write this using $k_i$ to den (meaning $u^p_k \notin N_{u^p_{k-1}}$).  % .Here $v$ cannot be a neighbor of $u$ or a neighbor of the last found special customer.

We now formally specify the set of LA-routes denoted $\Omega^{LA}$. Any $l \in \Omega^0$ lies in $\Omega^{LA}$ if the following holds for all $|l|\geq k_2>k_1\geq 1$ (with exposition below).%for which .  %is written as follows given any $k_2>k_1$.
\begin{align}
\label{ng_deg2}
    (u^l_{k_1}=u^l_{k_2}) \rightarrow \exists k_3\in Q^l, \mbox{s.t. } k_1<k_3<k_2, u^l_{k_1} \notin M_{u^l_{k_3}}
\end{align}
The premise of \eqref{ng_deg2} (left hand side of the $\rightarrow $ in \eqref{ng_deg2} ) is true if the same customer is at indexes $k_1$ and $k_2$.  The right hand side of \eqref{ng_deg2} states that there must exist a customer at special index $k_3$ that lies between $k_1$ and $k_2$ and does not consider $u^l_{k_1}$ to be a ng-neighbor. The set of LA-routes is a superset of the set of elementary routes and is a subset of the set of ng-routes.  The set of LA-routes is identical to the set of ng-routes when $|N_u|=0 \;  \forall u \in N$.   

% Later in the text we use $Q^l$ to denote the special indexes in any route $l \in \Omega$ or $l \in \Omega^{LA}$, which are defined using \eqref{qpdef}.  % and often index $l \in \Omega^{LA}$ similarly.
\end{itemize} 

We now consider an example of a ng-route that is not a LA-route. In our example, the set $N_u$ is identical to the set $M_u$ for each $u\in N$, $|N|=12$, and the locations of the twelve customers correspond to the set of positions on a classic analog clock. As a result, both the LA neighbors and the ng-neighbors of $u_k$ are [$u_{k-2}$,$u_{k-1}$,$u_{k+1}$,$u_{k+2}$] applied with modulus 12. Thus the neighbors (ng and LA) of $u_4$ (4 o'clock) are $N_{u_4}=[u_2,u_3,u_5,u_6]$ and the neighbors of $u_1$ are $N_{u_1}= [u_{11},u_{12},u_2,u_3]$. Respecting the properties set for our example, the following route is considered to be a feasible ng-route but is not a feasible LA-route: [-1,$u_3$,$u_1$,$u_5$,$u_1$,-2]. Observe that the set of special indices in this route consists of the index $1$, which corresponds to the customer $u_3$. 

\subsection{Local Area Structure of Reduced Cost}
\label{sec_pricing_graph}
We now describe costs and terms associated with computing the lowest reduced cost elementary route during pricing. % as done in \citep{mandal2022local}.
%Our description employs the following notation.  
%We use $N_u\subseteq N-u$ to denote the set of customers nearby $u$ for a given $u\in N$.  We refer to $N_u$ as the LA neighbors of customer $u$.  The selection of the cardinality of  $N_u$ involves computational trade offs seen later in the document though $|N_u|=10$ works well in practice. 
We use $\Omega_y$ where $y=(u,v,d)$ for some ($u\in N$, $v \in (N^+-(N_u\cup u \cup -1)),d_0-d_v\geq d\geq d_u$), which we index by $p$, to denote the set of elementary paths (which we call LA-arcs) meeting the following properties.  
    \textbf{(a):  } Path $p$ starts at $u$, ends at $v$ and all intermediate customers denoted $N_p$ lie in the LA neighborhood of $u$; meaning that  $N_p\subseteq N_u$.  
    \textbf{(b): }  The total demand serviced in path $p$ prior reaching $v$ is $d$ meaning $\sum_{w \in N_p\cup u}d_w=d$.
    \textbf{(c):  }Path  $p$ is the lowest cost path starting at $u$ ending at $v$ and visiting intermediate customers $N_p$.  We use $c_p$ to denote the total travel distance (cost) of path $p$. 
    %\textbf{(4):  }. 
    %The ordering of the intermediate customers $N_p$ in path $p$ is denoted $\zeta_p$ and is set so as to minimize $c_p$. 
We use $Y$, which we index by $y$, to denote the set of ($u\in N$, $v \in (N^+-(N_u\cup u \cup -1)),d\geq d_u$) for which $\Omega_y$ is non-empty. 
The LA-routes approach seeks to find a resource constrained shortest path on a directed acyclic multi-graph with the following properties. It is constructed so that any elementary path from source to sink has total cost on traversed edges equal to $\bar{c}_l$ and that the set of elementary paths from source to sink is exactly the set of elementary routes \footnote{Minus possibly some routes that have higher than necessary cost; meaning that they service the same customers as another route but have higher cost; hence such routes are not in any optimal ILP solution and can be ignored.}.  Thus finding the lowest reduced cost elementary path solves \eqref{pricing}.  We define this multi-graph as follows using $I,E$ to denote the set of nodes and edges respectively. 
%\begin{itemize}
%\item
We index $I$ using $i$ or $j$.  There is one node corresponding to the source $(-1,d_0)$ and the sink $(-2,0)$. There is one node for each $u,d$ where $u\in N$ and $d_0\geq d\geq d_u$. For any given edge $i,j \in E$ we define $N_{ij}$ to denote the set of customers serviced by traversing edge $i,j$.  % which contains all customers but the last customer in the Larc
%\item 
We connect $i=(-1,d_0)$ to $j=(u,d)$ for each $d\geq d_u,u\in N$.  This edge is associated with weight $\bar{c}_{ij}=c_{-1u}+\pi_0$. Traversing this edge indicates that the vehicle leaves the depot on a route servicing exactly $d$ units of demand, where its first customer is $u$.   Here $N_{ij}=\{ \}$.
 For each $y \in Y, p \in \Omega_y$ ( where $y=(u,v,d)$) and nodes in $i\in I,j\in I,$ s.t. $i=(u,d_1),j=(v,d_1-d)$ we connect $i$ to $j$ with an edge of weight $\bar{c}_{ij}=c_{p}-\sum_{w \in N_p \cup u}\pi_u$.  Traversing this edge indicates that upon arriving at $u$ with $d_1$ units of demand remaining that the vehicle travels on path $p$ servicing the customers in $N_p$ then proceeds immediately to $v$.  Here $N_{ij}=u\cup N_p$ (note that $v$ is not in $N_{ij}$).  
%\end{itemize}
%
Given a path from source to sink containing the arcs $\hat{E}$ this path is elementary (services no customer more than once) if the edges in $\hat{E}$ have disjoint $N_{ij}$ sets.  We use $a_{ijl}=1$ if route $l$ uses edge $ij$.  Thus we write the reduced cost of a route $l$ as follows in terms of $\bar{c}_{ij}$:  $\bar{c}_l=\sum_{ij \in E}a_{ijl}\bar{c}_{ij}$.  We use $p_{ij}$ to denote the path corresponding to edge $ij$ where $p_{ij}=0$ if $i=(-1,d_0)$.  
%\subsection{LA-routes:  Component for Pricing}
%\label{sec_la_routes_pricing}
%
%In this section we define LA-routes which we use to solve pricing.  We express LA-routes using the following terms  $M_u \subseteq N-u$  which we refer to as ng-neighbors as introduced by \citep{baldacci2011new}.  % and used by many different researchers
%The sets $M_u \forall u \in N$  grow over the course of DSSR \citep{righini2009decremental}. We use $\Omega^{LA}$ to denote the set of paths from source to sink in the graph in Section \ref{sec_pricing_graph} some of which are not elementary and hence $\Omega \subseteq \Omega^{LA}$.  To assist in our description we use $u^l_{k}$ to refer to the $k'th$ customer visited  in the route $l$. 
% which are the set of positions for customers in nodes visited on the LA-route in the multi-graph in Section \ref{sec_pricing_graph}.
%
\subsection{Decremental State Space Relaxation}
\label{sec_decremental_la}
Below we describe the use of Decremental State Space Relaxation (DSSR) to generate the lowest reduced cost elementary route.  %We define the $N_u$ for all $u\in N$ to be the composed of the nearest customers to customer $u$.
We initialize $M_u=\{ \}$  $\forall u \in N^+$. We then iterate between the following two steps until the lowest reduced cost LA-route generated is elementary.  
 \textbf{(1):  }Solve for the lowest reduced cost LA-route $l$, as a shortest path problem (not a RCSP) as to be described in Section \ref{sec_pricing_la}.
    \textbf{(2):  }
    Find a cycle of customers in route $l$ (if it exists). Consider that the cycle identified starts/ends with customer $u$ at indexes $k_1,k_2$ where $k_2>k_1$. 
    Now, for each  customer $k_3 \in Q^l$ for which $k_2> k_3>k_1$ add $u$ to $M_{u^l_{k_3}}$.  In Section \ref{sec_pricing_la} we shall see that smaller $M_u$ sets are computationally desirable for pricing over LA-routes (which is step 1 of DSSR).  Thus, intelligent cycle selection is done so as to keep $M_u$ sets small, which we discuss in Section \ref{astar_dssr}.    
In CG we need to only produce a negative reduced cost column at each iteration of pricing to ensure an optimal solution to the MP.
We map the non-elementary route generated at each iteration of DSSR to an elementary route and terminate DSSR when this elementary route has negative reduced cost.  In order to generate this elementary route, we remove each customer that is included more than once (after its first inclusion). Here we present the notation we will need for the next section. 
\begingroup
\setlength{\tabcolsep}{10pt} % Default value: 6pt
\renewcommand{\arraystretch}{1.5} % Default value: 1
\begin {table}[H]
\begin{center}
\begin{tabular}{|p{1.5cm}|p{3cm}||p{7cm}|}%|p{1cm}| p{2cm}|}
 \hline
Name & Space & Meaning \\%Comp Time & Performed in Pre-processing \\
 \hline
 $N_u$ & $N_u \subseteq \{N-u\}$ & The LA neighbors of $u$.  %These are used in the context of LA-routes.
 \\%& none & yes\\
 $M_u$ & $M_u \subseteq \{N-u\}$ & The ng-neighbors of $u$.% or decremental state space neighbors and inside LA-routes alongside $N_u$ 
 \\%& none & yes \\
 $z$ & $z \in Z$ & $z$ is a member of set $Z$ where $Z$ corresponds to the space of $u \in N^+, v \in N^+, M_1 \subseteq M_u, M_2 \subseteq M_v, d = d_1 - d_2$. %Each member $z$ lays out a path which obeys the following constraints:
 %The path starts at $u$ ends at $v$; visiting No customers in $M_1$; All customers in $M_2-M_1$ are visited where $M_2 \subseteq M_v$. No customers in $M_v-M_2$ are visited. Total demand serviced is  $d_1-d_2$  (excluding $v$) OR if $v$ is the sink then the total demand serviced does not exceed $d_1-d_2$. All customers visited lie in $N_u$ excluding $v$ which is the last customer in the path and does not lie in $N_u$.
 \\%& None & yes\\
 $i=(u,M_1,d)$& $u \in N, M_1 \subseteq M_u,d_0 - \sum_{w \in M_1}d_w \geq d\geq d_u$ &  Being at $i$ means that the nascent route is currently at $u$, has $d$ units of capacity remaining (prior to servicing $u$), and all customers in $M_1 \subseteq M_u$ have been visited and serviced at least once already. \\%.  & none & yes\\
  $\Omega_z$ & $z \in Z$&  The set of minimum cost elementary paths covering all possible sets of customers for elementary paths meeting the constraints set by $z=(u,v,M_1,M_2,d)$. \\%& fast & No (only updated when $M$ sets change given $\pi$)\\
 $p^z$ & $\Omega_z$ & The lowest reduced cost path in $\Omega_z$.\\% & fast & No (only updated when $M$ sets change or $\pi$ changes  \\
 $c_{p^z}$ & $p \in \Omega_z,z \in Z$ & The cost of the lowest reduced cost path in $\Omega_z$.\\% & Slow & yes.  Done once jointly for all $p$ Independent of $M$ sets\\
 \hline
\end{tabular}
\end{center}
\vspace{0.25cm}
\caption {\bf{Local Area Routes Notation}} %\label{tab_LANG} 
\label{tab_LANG}
\end{table}
\endgroup

\subsection{Pricing over LA-Routes}
\label{sec_pricing_la}
In this section we discuss generating the lowest reduced cost LA-route as required in step \textbf{(1)} in DSSR as a simple shortest path computation. To assist in this description we introduce the following terms. For any given tuple $z=(u,v,M_1\subseteq M_u,M_2\subseteq M_v,d)$ we use $\Omega_z$ to denote the subset of $\Omega_y$ where $y=(u,v,d)$ s.t. $p \in \Omega_z$ satisfies all of the following properties:  $|N_p\cap M_1|=0$, $M_2=M_v\cap (u \cup M_1 \cup N_p)$.  We use $Z$, which we index by $z$, to denote the set of non-empty $\Omega_z$. 
Given any dual solution $\pi$, we define $\bar{c}_z$ as the reduced cost of the lowest reduced cost LA-arc in $\Omega_z$ as follows (with minimizer denoted $p^z$):  
    $\bar{c}_z=\min_{p \in \Omega_z}\bar{c}_p$. 
%We define the minimizing LA-arc as $p^z$.  
Using terms $p^z$ we describe the following graph with vertex and edge sets denoted $I^M,E^M$ where the lowest reduced cost LA-route in $\Omega^{LA}$ corresponds to the lowest cost path from the source to the sink.
%\begin{itemize}
%
%\item
We index $I^M$ with $i,j$.  There is one node in $I^M$ corresponding to the source denoted $(-1,\{ \},d_0)$, and one node for the sink denoted $(-2,\{ \},0)$. 
For each node $(u,d) \in I$ (all nodes are defined in this manner in $I$ excluding source and sink), create one node $(u,M_1,d)$ in $I^M$ for each tuple $u,M_1,d$ satisfying the following: $ M_1 \subseteq M_u,d_0-\sum_{w \in M_1 }d_w\geq d\geq d_u$.
%
%\item
We connect $(-1,\{ \},d_0)$ to $(u,\{ \},d)$ for each $(u,\{ \},d)$ in $I^M$.  This edge is associated with weight $\bar{c}_{ij}=c_{uv}+\pi_0$.  Traversing this edge indicates that the vehicle leaves the depot on a route servicing exactly $d$ units of demand, where its first customer serviced is $u$. 
%\item
For each $z \in Z$ ( where $z=(u,v,M_1,M_2,d)$) and nodes in $i\in I^M,j\in I^M$ where $i=(u,M_1,d_1),j=(v,M_2,d_1-d)$ we connect  $i$ to $j$ with an edge with weight $\bar{c}_{ij}=c_{p^z}-\sum_{w \in N_{p^z} \cup u}\pi_u$.  Traversing this edge indicates that upon arriving at $u$ with $d_1$ units of demand remaining that the vehicle travels on path $p^z$ meaning that it services the customers in $N_{p^z}$ after leaving $u$ then proceeds immediately to $v$.  
%\end{itemize}

Given $I^M,E^M$ we can solve pricing using the Bellman-Ford algorithm to find the shortest path from -1 and -2 in $E^M$, since $E^M$ may have negative weights but no negative weight cycles as $E^M$ describes a directed acyclic graph.

The set of paths from source to sink in  $I^M,E^M$ does not include all LA-routes.  However we now establish that this set of paths includes the lowest reduced cost LA-route which is sufficient in order to complete step (1) of DSSR.  In order to permit every LA-route to be expressed we use a multi-graph $I^M,E^{M+}$ where $E^{M+}$ has all edges in $E^M$ plus the following additional edges.  For each  $z \in Z, p \in \Omega_z-p^z$ ( where $z=(u,v,M_1,M_2,d)$) and nodes in $i\in I^M,j\in I^M$ where $i=(u,M_1,d_1),j=(v,M_2,d_1-d)$ we connect  $i$ to $j$ with an edge with weight $\bar{c}_{ij}=c_{p}-\sum_{w \in N_p \cup u}\pi_u$.  Traversing this edge indicates that upon arriving at $u$ with $d_1$ units of demand remaining that the vehicle travels on path $p$ meaning that it services the customers in $N_p$ after leaving $u$ then proceeds immediately to $v$. 
Observe that the lowest cost path from source to sink in ($I^M,E^{M+}$), does not use any edge between a pair of nodes $i,j$ other than the one with lowest possible cost.  Thus the lowest  cost path in $(I^M,E^{M+})$ uses only edges in $E^M$. Therefore the lowest cost path on  $(I^M,E^{M})$ is the lowest reduced cost LA-route.  
\subsection{Fast Computation of LA-Arc Costs}
\label{subsubsec_LA_review}
In this subsection we describe the efficient computation of $c_p$ terms, which is done once prior to the first iteration of CG and is never repeated. % This section shows that $c_{p}$ terms are easy to compute using a dynamic program when the size of $N_u$ sets are small.
Let $P^+$ be defined as the set of tuples of the form $p=(u,v,\hat{N})$ where for each $p\in P^+$ there exists a  $u_1 \in N$ s.t. $(u \in (N_{u_1} \cup u_1),v \in N^+-(\hat{N}\cup u \cup u_1),\hat{N}\subseteq N_{u_1})$.  Let $c_{uv\hat{N}}$ for any $\{u,v,\hat{N}\} \in P^+$ denote the cost of the lowest cost path starting at customer $u$, ending at $v$, and visiting all customers in $\hat{N}$ which we denote as $c_{uv\hat{N}}$.  We write  $c_{uv\hat{N}}$ recursively below with helper terms $c_{u_2,u_2,\{\}}=0$ and $c_{u_2,v_2,\{\}}=c_{u_2,v_2}$ for all $u_2 \in N^+,v_2 \in N^+$ to describe the base cases. 
\begin{align}
    c_{u,v,\hat{N}}=\min_{w \in \hat{N}}c_{u,w}+c_{w,v,\hat{N}-w} \quad \forall (u,v,\hat{N}) \in P^{+}
    \label{my_sub_worker}
\end{align}
We generate the $c_{p}$ terms efficiently by iterating over the elements of $P^+$ in the order of increasing sizes of intermediate customer sets, using \eqref{my_sub_worker} to evaluate $c_{p}$. %
\subsection{ Accelerating DSSR with A*}
\label{astar_dssr}
In this section we exploit information from the first iteration of DSSR to solve subsequent iterations efficiently using the A* algorithm \citep{dechter1985generalized}. In order to apply A* we require the graph have non-negative weights. Let us offset $\bar{c}_{ij}$ where $i=(u,M_1,d_1)$ and $j=(v,M_2,d_2)$ by adding $\eta *(d_1-d_2)$ to $\bar{c}_{ij}$ where $\eta$ is the smallest value sufficient to ensure that all edge terms are non-negative.  Thus $\eta$ is defined as follows.  \begin{align}
\eta=-\min(0,\min_{\substack{p\in \Omega_y\\ y\in Y}}\frac{\bar{c}_p}{\sum_{w\in N_p \cup u}d_w})\end{align}
Observe that the cost of every path from source to sink is increased by exactly $\eta d_0$ so the minimizing cost path is not changed by this addition. 

The A* algorithm operates similarly to Dijkstra's algorithm except that it expands the un-expanded node $i\in I^M$ with the lowest sum of the distance to reach the node (denoted $g_i$) from the source and a heuristic ($h_i$). This heuristic $h_i$ is a lower bound on the cost of the lowest cost path to reach the sink starting from that node. We describe a high quality easily computed  heuristic as follows.  Given empty $M_u$ sets (as is the case during the first iteration of DSSR), we compute the shortest distance from each node $(u,\{ \},d)$  to the sink. We denote this heuristic as $h_{ud}$ and refer to the graph where all ng-neighbor sets are empty as the initial graph which we denote with nodes,edges  $(I^0,E^0)$.  We associate each $h_{ud}$ term to nodes of the form $i=(u,M_1,d)$ for each $M_1 \subseteq M_u$. The $h_{ud}$ terms are computed exactly (via Bellman-Ford) once for each call to pricing and not for each iteration of DSSR. This initial graph is much smaller than the graph in later iterations of DSSR so the computations of $h_{ud}$ is trivial.  
Observe that the $h_{ud}$ terms provide a consistent heuristic for A* since the set of LA-routes on the graph where all ng-neighbor sets are empty is a super-set of the set of LA-routes on any graph described by $(I^0,E^0)$ in DSSR when the ng-neighbor sets are non-empty.  

The A* search procedure motivates the following cycle selection criteria for DSSR.  In the step \textbf{(2)} of DSSR we select the cycle that causes the the least number of nodes to be added to the graph $(I^M,E^M)$, since adding extra nodes in the graph $I^M$ may require additional nodes to be expanded during A*. This is crucial as the number of nodes $(u,M_1,d)$ in the graph $(I^M,E^M)$ can grow exponentially in terms of $|M_u| \quad \forall u \in N$.
\subsection{Exploitation of Dominance Criteria}
\label{subsec_dom}
Dominance criteria are used in pricing algorithms to limit the number of expansion operations \citep{irnich2005shortest}.  To employ them we implement pricing as previously described over $(I^M,E^M)$ on a marginally different edge set than described thus far (with the same node set).   However the set of  LA-routes that can be represented in the original graph is the same as the LA-routes that can be represented in this modified graph and have identical costs.  We introduced the original graph formulation for pricing first as it is easier to understand and it is used in the stabilized CG algorithm discussed in Section \ref{sec_new_gg_pgm}. 
The modified graph changes the connections from the source and the connections to the sink only.  Both pricing and heuristic generation are done on this modified graph.  

 We connect the source only to nodes of the form $(u,\{ \},d_0)$ for each $u \in N$ with cost $\hat{c}_{ij}=c_{-1u}+\pi_u$ (for $i=(-1,\{ \},d_0)$,$j=(u,\{ \} ,d_0)$) in the new graph.  Traversing this edge indicates that the first customer visited on the route is $u$.  
 
 For each node $(u,M_1,d)$ the set possible LA-arcs associated with reaching $(-2,\{ \},0)$  includes all LA-arcs starting in $u$, visiting intermediate customers in $N_u-M_1$ with total demand (including $u$) not exceeding $d$, and ending at the depot. The cost of the edge is set so that every path from source to sink adds exactly $\eta d_0$ to $\bar{c}_l$.  %Edge costs are only modified on the sink connections, and for all other connections, the edge costs match the original graph. Furthermore, connections from a node $u,M_1,d$ to the sink now must consider all LA-arcs using up to $d$ units of demand; and add the the existing cost along the edge to $\eta$ times the capacity remaining at the end of the LA-route (so as to ensure that all routes have added cost of exactly $\eta d_0$).  %The preiovus presenmta%The modification easily permits the use of dominance criteria accelerating pricing \citep{irnich2005shortest}.  In this graph all nodes are connected to the sink and fewer are connected to the source. 
Thus we define the updated $\hat{c}_{ij}$ for sink connections as follows.  
\begin{align}
\label{get_new_cost}
\hat{c}_{ij}\leftarrow d\eta+ \min_{\substack{d\geq d_1\geq d_u\\p \in \Omega_z\\z=(u,-2,M_1,\{ \},d_1)}}\bar{c}_p, \quad \forall i=(u,M_1,d),j=(-2,\{ \},0)
\end{align} 
The LA-arc associated with traversing the edge between $(ij)$ where $i=(u,M_1,d),j=(-2,\{ \},0)$ becomes the arg minimizer of \eqref{get_new_cost} and is denoted $p_{ij}$.  Traversing this arc indicates that the vehicle leaves $u$ visits the intermediate customers of $p_{ij}$ then returns to the depot.  %Observe that 
%The assocaited arc is$
 %Computing the lowest reduced cost LA-route on this graph proceeds as in the original form; however we can now make the following observation.
 
 Consider any $u,M_1\subseteq M_2\subseteq M_u,d_u\leq d_2\leq d_1$.  Observe that the lowest reduced cost path staring at $(u,M_1,d_1)$ and ending at the sink can have reduced cost no greater than the lowest reduced cost path starting at $(u,M_2,d_2)$ and ending at the sink.  This is because any such sequence of customers  starting at $(u,M_2,d_2)$ is also valid for $(u,M_1,d_1)$.  Thus we need never expand a node $(u,M_2,d_2)$ in $A^*$ if we have already expanded the node $(u,M_1,d_1)$.  %the distance from the source to $(u,M_1,d_1)$ is less than the distance from the source to $(u,M_2,d_2)$.

\section{LA-Arcs Encoding SRI}
\label{Sec_SRI_LA}
In this section we relax subset row inequalities (SRI) \citep{jepsen2008subset} to produce LA-SRI so as to permit efficient inclusion of SRI (in a marginally weakened form) in our LP relaxation and pricing.  
For any $y\in Y,p \in \Omega_y,\delta \in \Delta$ where $y=(u,v,d)$ let us define $a_{p\delta}$ as follows:  $a_{p\delta}=\sum_{w \in N_p \cup u}[w \in N_{\delta}]$.  Let $a_{pl}=1$ if LA-arc $p$ is used in route $l$. We now write SRI in terms of $a_{p\delta},a_{pl}$ then take a lower bound to produce LA-SRI.  
\begin{align}
\label{lasrieq}
\sum_{l\in \Omega}\theta_l 
\left \lfloor\frac{\sum_{\substack{y\in Y\\p\in \Omega_y}}a_{p\delta}a_{pl}}{m_{\delta}}\right \rfloor
\geq \sum_{l\in \Omega}\sum_{\substack{y\in Y\\p\in \Omega_y}}\theta_l a_{pl} \left \lfloor\frac{a_{p\delta}}{m_{\delta}}\right \rfloor
\quad \forall \delta \in \Delta 
\end{align} 
 Below we write LA-SRI for a given $\delta$ (with associated dual variable $\pi_{\delta}$) using helper terms $a^*_{p\delta},a_{\delta l}^*$ defined  as $a^*_{p\delta}=\left \lfloor\frac{a_{p\delta}}{m_{\delta}}\right \rfloor $; and $
     a^*_{\delta l}=\sum_{\substack{p \in \Omega_y \\ y \in Y}}a_{pl}a^*_{p\delta}$.% \quad \forall l \in \Omega, \delta \in \Delta$:
 %\begin{subequations}%\quad \forall p \in \cup_{y \in Y}\Omega_y, \quad \delta \in \Delta
 \begin{align}
 \sum_{l \in \Omega}a^*_{\delta l}\theta_l\leq \left \lfloor \frac{|N_{\delta}|}{m_{\delta}}\right \rfloor \quad \forall \delta \in \Delta \quad  [-\pi_{\delta}] \label{LASRIEQ}%\label{CVRP_valid_2}
 \end{align}
 % \end{subequations}
Observe that \eqref{LASRIEQ} can be understood as replacing constraints on routes as defined in equation \eqref{my_div_2} with constraints on LA-arcs. For example, for a LA-SRI $\delta$ with $m_{\delta}=2,|N_{\delta}|=3$, enforces that the number of LA-arcs including 2 or more customers in $N_{\delta}$ in the solution can not exceed one.  %$Observe that $a^*_{\delta l}$.  
We can now solve pricing as is done in in Section \ref{sec_LA} by incorporating $\pi_{\delta}$ into the $\bar{c}_p$ terms as follows:  $\bar{c}_p=c_p-\sum_{w \in N_p\cup u}\pi_w+\sum_{\delta \in \Delta_R}\pi_{\delta} a^*_{p\delta}$.  

 We now contrast LA-SRI with SRI in two examples, with further details in Appendix \ref{sec_diff_LASRI_sri}. In our examples, we use the case of a depot in San Diego (SD) and customers in NYC as defined for SRI except $3$ additional customers of unit demand are in SD denoted SD1, SD2, SD3 which are nearly co-located. 
\begin{itemize} 
\item \textbf{Example One: } LA-SRI do not constrain as well as SRI when customers in the SRI are not localized in space.  Consider a case with SRI defined to be $N_{\delta}$=[NYC1,NYC2,SD1] and $m_{\delta}=2$. Now consider routes $l_1$ and $l_2$ having sequences [-1, NYC1, NYC3,  SD1, -2]  and [-1, NYC1, NYC2, SD1, -2] respectively.  Observe that $a_{\delta l_1}=1$ and $a^*_{\delta l_1}=0$ and $a_{\delta l_2}=1$ and $a^*_{\delta l_2}=1$. \\   
\item \textbf{Example Two:  } LA-SRI do not constrain as well as SRI when a route returns to an area of localized customers after visiting it. However it is important to note that an optimal solution is disinclined to use such routes as they would tend to have prohibitively high cost. Consider a case with an SRI $\delta$ defined as  $N_{\delta}$=[NYC1, NYC2, NYC3] and $m_{\delta}=2$. Now consider routes $l_3$ and $l_4$, which have sequences [-1, NYC1, NYC2, SD1, -2] and, [-1, NYC1, SD1, NYC2, -2] respectively. Note that the special indexes are 1,3 for $l_3$ and 1,2,3 for $l_4$. Observe that $a_{\delta l_3}=a^*_{\delta l_3}=1$ and $a_{\delta l_4}=1$ and $a^*_{\delta l_4}=0$. 
\end{itemize}

%Additionally, a key opportunity for the LA-SRI to be weaker than the corresponding SRI, when the customers of the SRI is localized in space, is when LA-arcs are used on the route for which the following holds: the penultimate customer of the LA-arc is nearby the final customer of the LA-arc.  

%Additionally, a situation in which the LA-SRI could be weaker than the corresponding SRI is when the customers of the SRI are localized in space and when the LA-arcs are used on a route for which the penultimate customer of the LA-arc is nearby the final customer of the LA-arc.

Additionally, a situation in which an LA-SRI (with customers localized in space) could be weaker than the corresponding SRI is when an LA-arc is used for which the penultimate customer of the LA-arc is nearby the final customer of the LA-arc.

%Another possibility is when routes 
%

\section{Stabilization Exploiting LA Structure}
\label{sec_new_gg_pgm}
In this section we describe a stabilized version of CG adapted to efficiently solve the CG master problem (MP) that does not alter pricing. We achieve this by altering the sequence of dual solutions generated by CG algorithm with the aim that fewer iterations of CG are required to solve the MP. At each iteration of CG a more computationally intensive RMP is solved that includes a larger set of columns than those generated during pricing. Solving over this set of columns does not dramatically increase the time taken to solve the RMP since a simple computational structure is used to ensure that an exponential number of such columns can be encoded with a finite number of variables.  This set of columns is written as $\cup_{l \in \Omega_R}\Omega_l$ where $\Omega_l$ contains an exponential number of columns related to column $l$.  We refer to $\Omega_l$ as the family of column $l$. 
As a result, when solving the RMP in our stabilized version of CG, we solve $ \Psi(\cup_{l\in \Omega_R}\Omega_l,\Delta_R)$ instead of $\Psi(\Omega_R, \Delta_R)$ to accelerate the convergence of CG.  In this section we develop an efficient solver for this approach that does not enumerate all $\cup_{l\in \Omega_R}\Omega_l$ and does not employ expensive pricing operations.   
%We organize this section as follows.  In Section \ref{sub_sec_family} we formally characterize the families in our formulation and  the generation of a family from a column.  In Section \eqref{subsec_new_RMP} we describe the new RMP over graphs associated with families of columns leading to faster convergence inspired by CG. 
%

%\label{sub_sec_family}
%
\textbf{Families of Columns for CVRP: }For any given route $l$ we associate it with a strict total order of the members $N^+$ denoted $\beta^l$ where $-1,-2$ have the smallest and largest $\beta^l$ values respectively.  A route $\hat{l} \in \Omega$ lies $\Omega_l$ if and only if LA-arcs consistent with the ordering $\beta^l$ are used.  
An LA-arc $p$ in $\Omega_{y}$ (for y=$(u,v,d)$) is consistent with the ordering $\beta^l$ if $\beta^l_u < \beta^l_w<\beta^l_v$ for any $w \in (N_p)$; and $\beta^l_u<\beta^l_v$.  It is mandatory that route $l$ lies in $\Omega_l$ for all $l \in \Omega$.  
 We use $\Omega^l_y$ to refer to the subset of $p \in \Omega_y$ for which $p$ is consistent with ordering $\beta^l$.  We use $Y^l$ to denote the subset of $Y$ (indexed by $y$) for which $|\Omega^l_y|\neq 0$.  We use $E^l$ to denote the edges connected either to the source or for which the corresponding $y$ has non empty $\Omega^l_y$.  
 
 We now describe a procedure to generate $\beta^l$  motivated by the observation that customers that are in similar physical locations should be in similar positions on the ordered list. Having customers in such an order permits routes in $\Omega_{l}$ to be able visit customers close together without leaving the area and then coming back. We use a set $N_l$ to denote the customers in the route $l$. We initialize the ordering with the customers in $N_l$ sorted in order from first visited to last visited. Then, we iterate over $u\in (N-N_l)$ (in a random order) and insert $u$ behind the customer in $N_l$ nearest to $u$.  We insert $u$ that are closer to the depot than any customer in $N_l$ at the beginning of the list.  Observe that by using the aforementioned construction that $l$ lies in $\Omega_l$ for all $l \in \Omega$. We use $E^l \subseteq E^0$ (where node/edge sets $I^0,E^0$ is defined in Section \ref{astar_dssr}) to denote the set of edges in $E^0$, $(u,d_1),(v,d_2)$ for which $\beta^l_u<\beta^l_v$ and in addition where either $u=-1$ or $\Omega^l_{y=(u,v,d_1-d_2)}$ is non-empty.
 
\textbf{Stabilization of RMP:  }
%subsec_new_RMP\label{subsec_new_RMP}
We now describe an efficient LP for $\Psi(\cup_{l \in \Omega_R} \Omega_l, \Delta_R)$.  We now define an LP over the graphs defined by the nodes,edge sets $I^0,E^l$ for each $l \in \Omega_R$. 
We set decision variable $x^l_p=1$ to use LA-arc $p$ for creating a route in family $\Omega_l$.  % be a variable which denotes if a $p$ composing a path in the family $l$ is used in the solution, and how many vehicles (can be fractional or 1) are used for $p$.
There is one such variable $x^l_p$ for each $p \in \cup_{y \in Y^l}\Omega^l_y,l \in \Omega_R$. % for some  $y \in Y^l$.
We set $x^l_{ij}=1$ to select edge $ij$ for creating a route in family $\Omega_l$.  There is one such variable for each $l\in \Omega_R,ij \in E^l$.  Using our decision variables $x,\theta$ we write $\Psi(\cup_{l \in \Omega_R}\Omega_l,\Delta_R)$ as a LP below which we refer to as $\Psi^+(\Omega_R,\Delta_R)$ (with $\Psi(\cup_{l \in \Omega_R}\Omega_l,\Delta_R)=\Psi^+(\Omega_R,\Delta_R)$ proven in Appendix \ref{subseqEq}). % \ref{subseqEq}).%  We provide annotation of this RMP below the equations. 
\begin{subequations}
\label{myRMPSSTAB}
\begin{align}
    \min_{\substack{x\geq  0\\\theta \geq 0}}\sum_{l \in \Omega_R}c_l\theta_l +
    \sum_{\substack{l \in \Omega_R\\y \in Y^l\\ p \in \Omega^l_y}}c_{p}x^l_p+\sum_{\substack{l \in \Omega_R\\ u \in N\\ d_0\geq d\geq d_u}}c_{-1 u}x^l_{(-1,d_0),(u,d)}
    \label{RMP_stab_obj}\\
    \sum_{l \in \Omega_R}\theta_l+
    \sum_{l \in \Omega_R}\sum_{ij \in E^l}[i=(-1,d_0)]x^l_{ij}\leq K\quad [-\pi_0] \label{RMP_stab_packing}\\
\sum_{l \in \Omega_R}a_{ul}\theta_l+
\sum_{\substack{l \in \Omega_R\\y \in Y^l\\ p \in \Omega^l_y}}x^l_{p}a_{up}\geq 1 \quad \forall u \in N\quad [\pi_u] \label{RMP_stab_cover}\\
    \sum_{l \in \Omega_R}a^*_{\delta l}\theta_l+
    \sum_{\substack{l \in \Omega_R\\y \in Y^l\\ p \in \Omega^l_y}}a^*_{\delta p}x^l_p\leq \left \lfloor \frac{|N_{\delta}|}{m_\delta}\right \rfloor \quad \forall \delta \in \Delta_R \quad [-\pi_{\delta}]\label{RMP_stab_delta}\\
    \sum_{p \in \Omega^l_y}x^l_{p}= \sum_{\substack{ij\in E^l\\i=(u,d_1)\\ j=(v,d_1-d)}}x^l_{ij}\quad  \label{RMP_stab_agree} \forall y \in Y^l; y=(u,v,d),l \in \Omega_R\\
    \sum_{ij \in E^l}x^l_{ij}=\sum_{ji \in E^l}x^l_{ji} \quad \forall i \in I^0,  u_i\notin \{-1,-2\},l \in \Omega_R \label{RMP_stab_flow}
\end{align}
\end{subequations}
    \eqref{RMP_stab_obj}:  We seek to minimize the sum of the costs of the LA-arcs and edges (from the depot) and routes used in our solution. 
     \eqref{RMP_stab_packing}:  We enforce that no more that $K$ vehicles leave the depot.
     \eqref{RMP_stab_cover}:  We enforce that each customer is serviced in at least one LA-arc or route.
     \eqref{RMP_stab_delta}:  We enforce all LA-SRI over the LA-arcs and routes which form our solution. 
      %This means that for every family $l$, that the total weight of incoming edges is equal to that of the outgoing edges on all nodes corresponding to customers $u \in N$.
      \eqref{RMP_stab_agree}:  We enforce that the the LA-arcs selected are consistent with the edges selected for a given family.
      \eqref{RMP_stab_flow}:  We enforce that the edges selected for each family $l$ describe a set of routes using $u_i$ to denote the customer or depot associated with $i$.
      
The CG optimization procedure used to solve $\Psi(\Omega,\Delta_R)$ alternates between solving \eqref{myRMPSSTAB} (meaning $\Psi^+(\Omega_R,\Delta_R)$) and adding to $\Omega_R$ the lowest reduced cost column $l \in \Omega$.  We price over $\theta$ (never $x$) to add in a new column $l$ to $\Omega_R$. Then, we add to the RMP the $x$ terms associated with $l$ (which are $x^l_{ij}$ for all $ij \in E^l$ and $x^l_p$ for all $p\in \Omega^l_y,y \in Y^l$). We never consider $\theta$ terms when solving the RMP as they are redundant given $x$.  However we include the $\theta$ terms in \eqref{myRMPSSTAB} so as to make clear that pricing for \eqref{myRMPSSTAB} is identical to previous pricing problems. %  pricing operation being identical to teh previous RMP

Our solver for $\Psi(\Omega,\Delta)$ applies CG and cutting plane approaches (cut-price). Specifically  we solve for $\Psi(\Omega,\Delta)$ in a cutting plane manner where we add after each solution to $\Psi(\Omega,\Delta_R)$ the most violated LA-SRI (in experiments we chose the $30$ most violated LA-SRI).  %, or all LA-SRI when there are less than $30$).
This requires that we solve for $\Psi(\Omega,\Delta_R)$ which we solve using CG with the RMP being $\Psi^+(\Omega_R,\Delta_R)$.  % 
%Given our solver for $\Psi^+(\Omega_R,\Delta_R)$  we solve for $\Psi(\Omega,\Delta_R)$  using CG with $\Psi^+(\Omega_R,\Delta_R)$ as the RMP and pricing solved for using the solver in Section \ref{sec_LA}.  Given our solver for $\Psi(\Omega,\Delta_R)$
%
%
Observe that  $\Psi^+(\Omega_R,\Delta_R)$ may have a very large number of variables as $|\Omega_R|$ grows and for problems where $|\Omega_y|$ is large for some $y \in Y$; thus making the solution to \eqref{myRMPSSTAB} at each iteration of CG intractable.  %Hence we solve \eqref{myRMPSSTAB} in a manner akin to CG.
Thus we seek to construct a small set of edges, LA-arcs  denoted $\hat{E}^l\subseteq E^l \quad \forall l \in \Omega_R,\quad \hat{\Omega}_{y}^l \subseteq \Omega_{y}^l$ \quad $(y \in Y^l \; \forall l \in \Omega_R)$ respectively s.t. solving \eqref{myRMPSSTAB} over these terms (denoted $\Psi^+(\Omega_R,\Delta_R,\{ \hat{E}^l\},\{\hat{\Omega}^l_y\}))$ yields the same solution as \eqref{myRMPSSTAB};  where for short hand we use $\{ \hat{E}^l \}, \{ \hat{\Omega}^l\} $ to describe the $\hat{E}^l,\hat{\Omega}^l_y$ for each $l \in \Omega_R,(y\in Y^l,l \in \Omega_R)$.  High quality integer solutions can be obtained by efficiently by solving $\Psi^+(\Omega_R,\Delta_R,\{ \hat{E}^l \}, \{ \hat{\Omega}^l_y \} )$ as an ILP at termination of CG.  An exact solution can be pursued by incorporating $\Psi(\Omega,\Delta)$ into a branch-price \citep{barnprice} solver.
 To solve \eqref{myRMPSSTAB} we alternate between the following two steps (details in the Appendix \ref{subsec_eff_solve}) which we refer to as  RMP-LP and RMP-Shortest Paths respectively.  
\textbf{Step One:}   Solve small LP $\Psi^+(\Omega_R,\Delta_R,\{ \hat{E}^l \}, \{ \hat{\Omega}^l_y \} )$ producing the dual solution  $\pi$.
    \textbf{Step Two:}  Iterate over $l \in \Omega_R$ and compute the lowest reduced cost route in $\Omega_l$ denoted $\hat{l}^l$ as a simple shortest path computation (not a RCSP).  The associated graph is  $(I^0,E^l)$ with edge weights defined in Section \ref{sec_pricing_la} except only considering $\Omega^l_y$ not $\Omega_y$ meaning  $\bar{c}_{ij}=\min_{\substack{p \in \Omega^l_y\\ y=(u,v,d_1-d_2)}}\bar{c}_p $.  
    We denote the edges/LA-arcs used in $\hat{l}^l$ as  $E^l(\hat{l}^l)$ and  $ \Omega^l_y(\hat{l}^l) $ for each $y \in Y^l$ respectively.  When $\bar{c}_{\hat{l}^l}<0$ we augment each $\hat{E}^l$ with  $E^l(\hat{l}^l)$ and augment each $\hat{\Omega}^l_y$ with $\Omega^l_y(\hat{l}^l)$  $(\forall y \in Y^l)$.  When $\bar{c}_{\hat{l}^l}\geq 0 \;  \forall l \in \Omega_R$ we have solved \eqref{myRMPSSTAB} optimally and terminate optimization.
    %\end{itemize}
%
 %
\section{Experiments}
\label{sec_exper}
In this section we demonstrate the value of LA-routes and our stabilization approach for the efficient solution and tightening of the ELP on CVRP.%  We consider a dataset of CVRP problem instances that vary by numbers of customers and vehicle capacity. We perform experiments using our approach from Section \ref{subsec_new_RMP} with different parameters for the LA-SRI (meaning the of range of $|N_{\delta}|$ terms for $\delta \in \Delta$) and the sizes of LA neighbor sets (sizes of $N_u$). %Larger sizes of sets  $\Delta$,$|N_u|$ have the potential to tighten the expanded LP relaxation, but increase the amount of computation time.%  In these experiments, pricing and determination of the lowest reduced cost route is solved using LA-routes and DSSR as described in Section \ref{sec_decremental_la}. \citep{righini2008new,righini2009decremental}.  
%We organize this section as follows. In Section \ref{exper_sub_alg} we consider the parameterizations used for our CG solver and describe our CVRP problems in our data set.
%In Section \ref{Experiments} we conduct experiments to demonstrate the value our approach for tightening and efficiently solving the expanded LP relaxation.
%\subsection{Algorithms and Data Sets Compared}
%\label{exper_sub_alg}
%
\subsection{Synthetic Data}
We now describe our parameterizations for problems in our dataset and for the CG solver. We considered CVRP problem instances with customers all having unit demand, and the following parameter possibilities for (number of customers,vehicle capacity): (20,4);(40,8). We generate 10 problem instances for each CVRP parameter possibility. Each instance randomizes the locations of customers, the starting depot, and ending depot.  We used a total of twelve parameterizations for our CG solver. Larger LA neighbor sets or larger sets in $\Delta$ (as defined by LA-SRI parameters) can tighten the LP relaxation but may increase computation time. Each CG parameterization is associated with a number of LA neighbors per customer, and the type of LA-SRI used. We set LA neighbor sets for a given customer $u$ with size $x$ to be the $x$ closest customers excluding $u$ (by spatial distance) to $u$. We considered four choices for the number of LA neighbors per customer ($4,6,8,10$) and three choices for LA-SRI parameter values denoted $a,b,c$ which are described below using: $\Delta_{\alpha,\beta}$ to denote the set of LA-SRI consisting of $\alpha$ (meaning $|N_{\delta}|=\alpha$) unique customers in $N$ and for which  $m_{\delta}=\beta$.  
\begin{itemize}
    \item Option a: $\Delta\leftarrow \Delta_{3,2}$. 
    \item Option b: $\Delta\leftarrow \Delta_{3,2} \cup \Delta_{4,3}$. 
    \item Option c: $\Delta\leftarrow \Delta_{3,2} \cup \Delta_{5,2} \cup \Delta_{5,3}$.
\end{itemize}
Note that Option c does not include $\Delta_{4,3}$ because $\Delta_{5,3}$ dominates  $\Delta_{4,3}$. 

We quantify the amount of tightening of the LP relaxation using a measure of ``relative increase", which describes the proportion of the gap between the optimal integer solution computed (over parameterizations defining the CVRP problem(s) and the CG solver) and the lower bound (the LP value) that is closed by adding LA-SRI. Thus a relative increase of 1 indicates that the CG solver closes all of the gap while a relative increase of 0 indicates no improvement at all. All experiments were run on a 2014 Macbook pro running Matlab 2016.
%\subsection{Experiments Conducted}
%\label{Experiments}
%

\begin{figure*}[!hbtp]
\begin{minipage}{1\textwidth}
\includegraphics[width=0.5\linewidth]{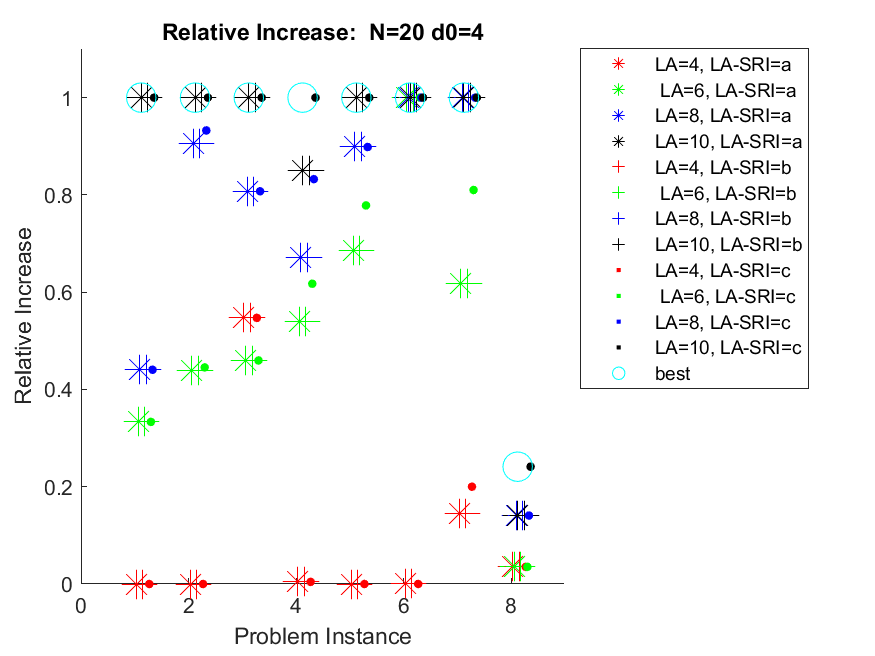}
\includegraphics[width=0.5\linewidth]{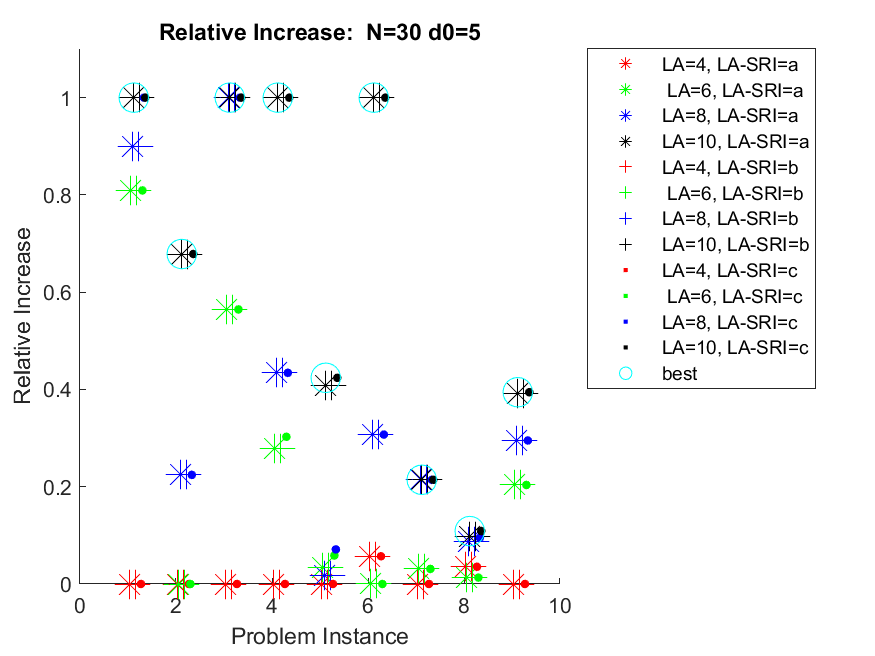}
\includegraphics[width=0.5\linewidth]{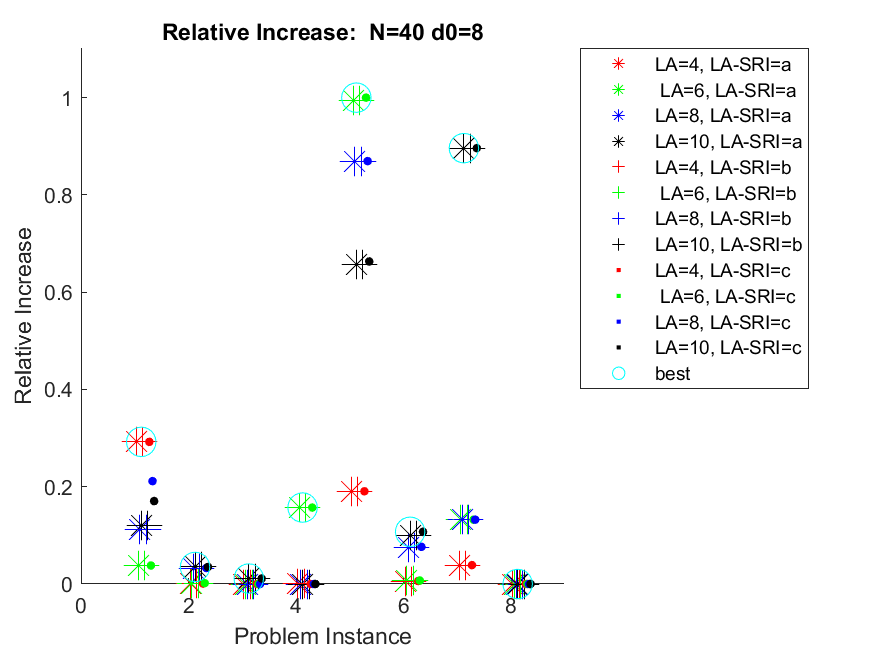}

\caption
{
The relative solution tightness for each CG solver/problem parameterization.
} 
    
\label{myCompTime}
 \end{minipage}
\end{figure*}

Figure \ref{myCompTime} shows the relative solution tightness for each CG solver/problem parameterization.% The y axis describes the relative increase in the tightness by a given solver parameterization. 
 Each data point describes the ``relative increase" in tightness for the CG solver for the problem instance. There is also a blue circle to indicate the greatest relative increase for a particular problem instance. For these problems only we ignore instances where the CG solver provides a tight solution before introducing LA-SRI, as the relative increase would be undefined. 

\newpage
 
\begin{figure*}[!hbtp]
\begin{minipage}{1\textwidth}

\includegraphics[width=0.5\linewidth]{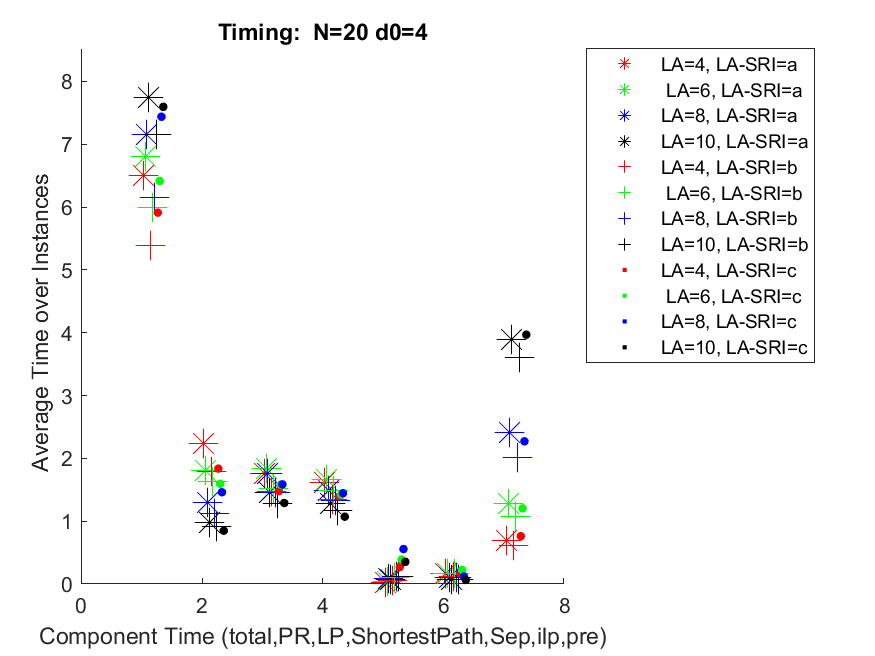}
\includegraphics[width=0.5\linewidth]{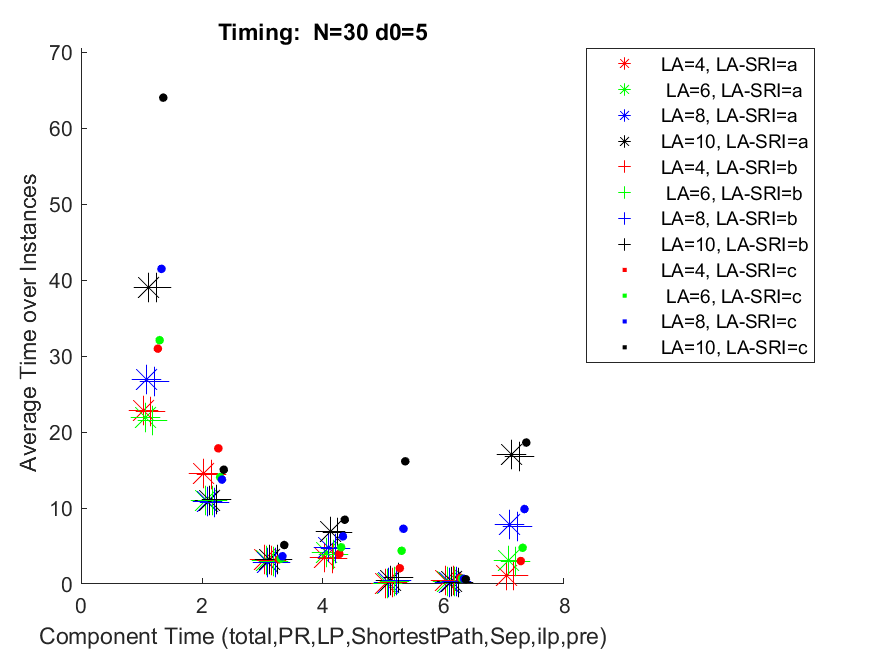}
\includegraphics[width=0.5\linewidth]{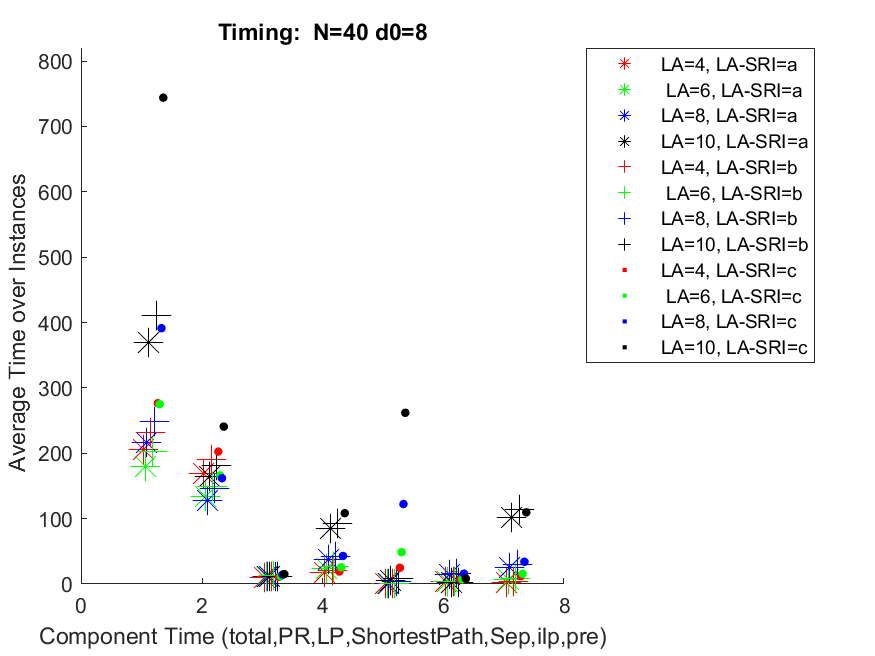}

\caption
{
The average time taken by a given CG solver parameterization for a given component of optimization.
} 
    
\label{myCompTime2}
 \end{minipage}
\end{figure*}

 In Figure \ref{myCompTime2} each data point describes the average time taken by a given CG solver parameterization for a given component of optimization. The components are along the X axis and are shown in Table \ref{table:Components}.

\begin {table}%[H]
\begin{center}
\begin{tabular}{ |c|c| } 
 \hline
 x-axis & Component\\ 
 \hline
 1 & Total\\ 
 2 & Pricing\\ 
 3 & RMP-LP\\ 
 4 & RMP-Shortest Path\\ 
 5 & Separation of SRI\\ 
 6 & Solving ILP\\ 
 7 & Preprocessing to generate the $c_p$ terms  \hspace{2cm}\\
 \hline
\end{tabular}
\end{center}
\vspace{0.20cm}
\caption{Algorithm Components}
\label{table:Components}
\end{table}

%
%\caption {\bf{CVRP Route Notation}} 
%\label{basicTab}
%
% 
%
%%\;\forall i=(u,d_1),j=(v,d_2), ij \in E^l$.
    %Observe that via the ordering of $\beta^l$ that the route $\hat{l}^l$ must be elementary as all paths from source to sink describe elementary routes.  
    %If $\hat{l}^l$ has negative reduced cost then we update the edges, LA-arcs associated with route $l$
%
%\newpage

\begin{figure*}[!hbtp]
\begin{minipage}{1\textwidth}

\includegraphics[width=0.5\linewidth]{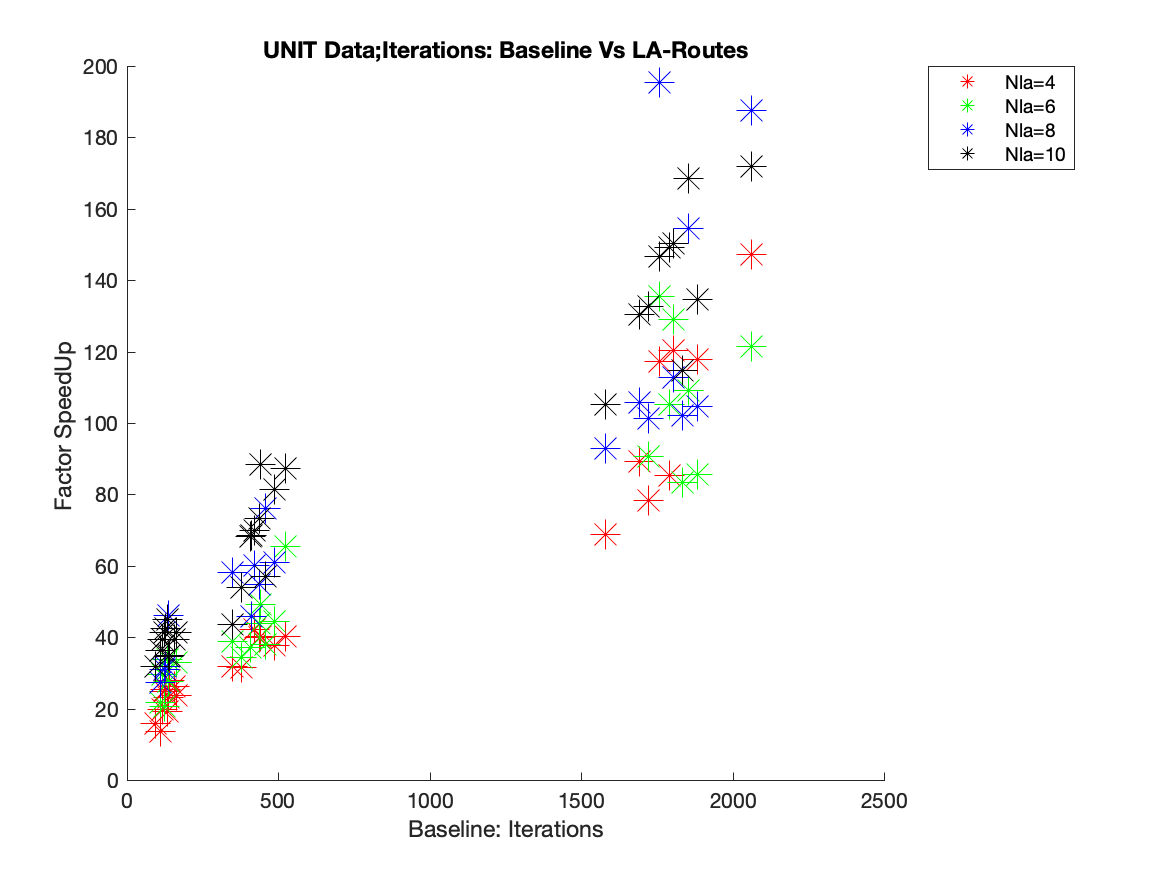}
\includegraphics[width=0.5\linewidth]{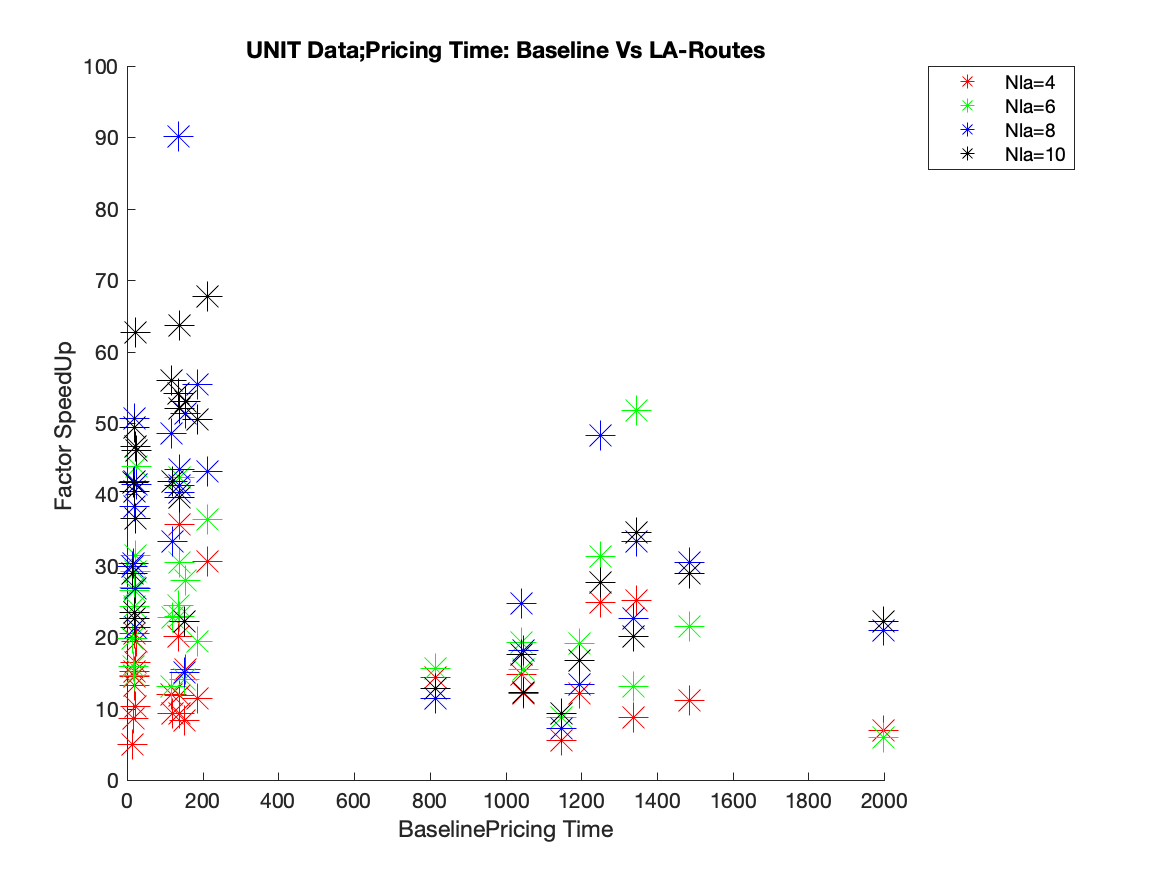}
\includegraphics[width=0.5\linewidth]{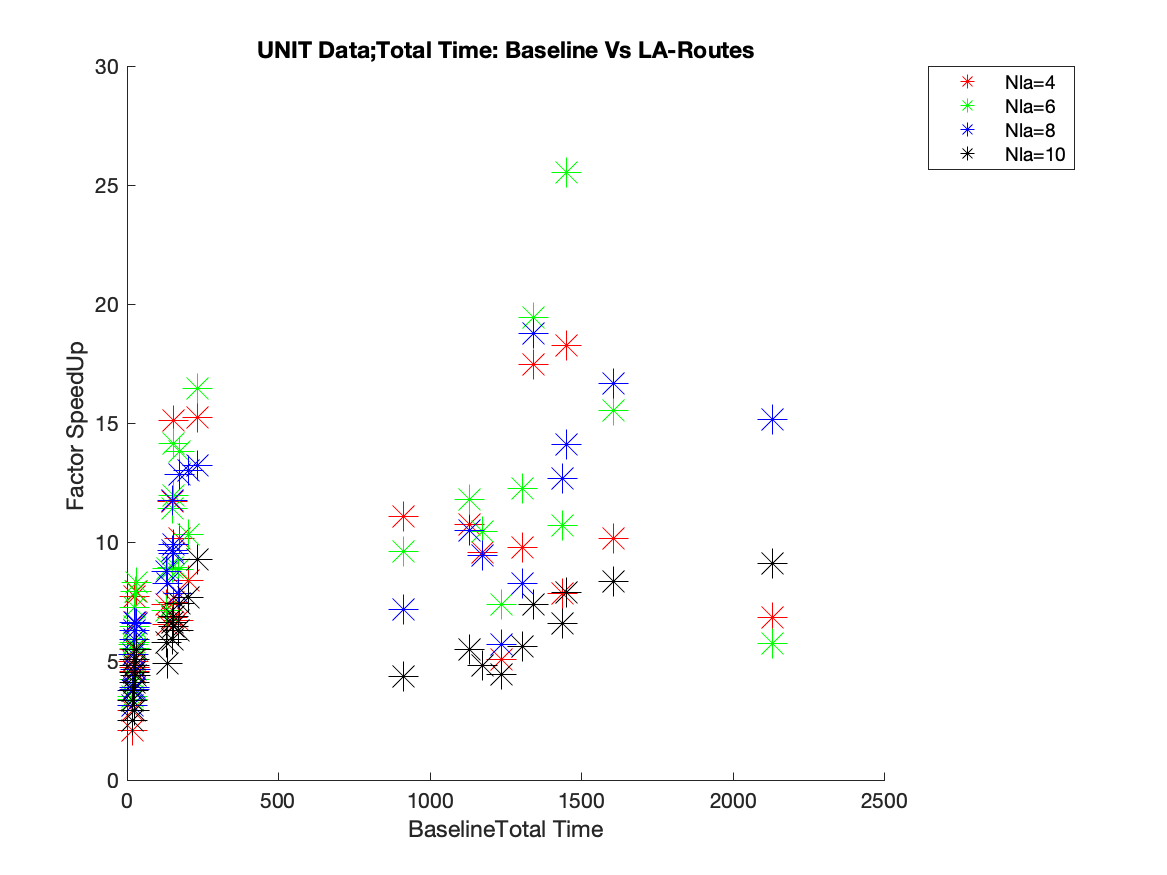}

\caption
{
The speed up relative to the baseline in terms of iterations, pricing time and total time.
} 
    
\label{myCompTime3}
 \end{minipage}
\end{figure*}

The third category, represented by plots in Figure \ref{myCompTime3}, contrasts the time taken to solve the MP without LA-SRI relative to a benchmark solver using DSSR for pricing and no stabilization. 

We observe large speed ups in terms of iterations, pricing time and total time using our approach over the baseline. We explore additional synthetic data with non-unit demands in Section \ref{secAddExper}.  %We also observe that pre-processing is a large consumer of computation time in some settings.  
\subsection{Benchmark Instances from Augerat A}
Here we examine our method on the problems referred to as Augerat A \citep{augerat1995computational}. We replicate the plots in Fig \ref{myCompTime3} for data set Augerat A. Note that in order to solve these problems quickly we make a small modification in which we divide the demands (and vehicle capacity) by five or ten (separate experiments are performed for each possibility) and then set the new demand (and capacity) to the ceiling of the values obtained. Since the baseline had trouble solving these problems we capped computation time for the baseline approach at 3000 seconds. Such instances treat the baseline as if it had finished solving the master problem. No time cap is used for our approach. We observe large speed ups in terms of iterations, pricing time and total time.  

In Appendix \ref{RCI_Augerat} we examine the performance of our approach on Augerat A, with rounded capacity inequalities (RCI)\citep{archetti2011column} and LA-SRI; along with a brief discussion of RCI in the context of LA-arcs.  %is . 
%\par To understand the figure imagine that the performance of each of the test cases are shown in vertical lines containing the four identifiers for LA-routes of sizes 4, 6, 8 and 10. In most cases the number of iterations will be lower with larger LA-routes and the number of iterations will also be lower. However the speed-ups in overall time are lower proportionally because of the overhead involved in implementing the pre-processing which can be sped up and parallelized in an industrial version.  %method generally and then in the exponential increase in terms pre-processing (generation of $c_{p}$ terms) of the computation that the larger LA-routes cause; which can be efficiently sped up and parallelized. Therefore, while we know that we can speed up computation by 8-10X, we are confident that we can improve on these results. In addition we believe that the computation time with larger LA-routes will be reduced relative to those with smaller LA-routes.  
\newpage
\begin{figure*}[!hbtp]
\begin{minipage}{1\textwidth}
\includegraphics[width=0.5\linewidth]{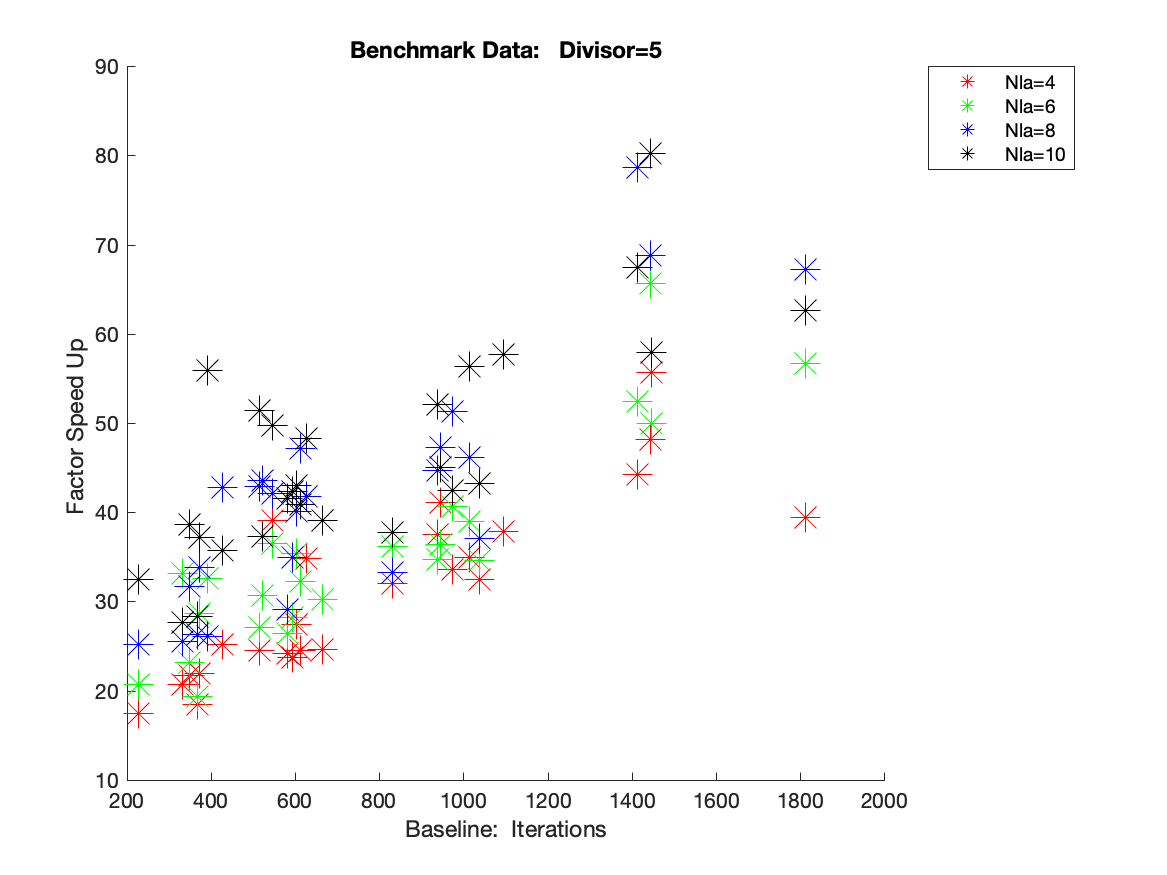}
\includegraphics[width=0.5\linewidth]{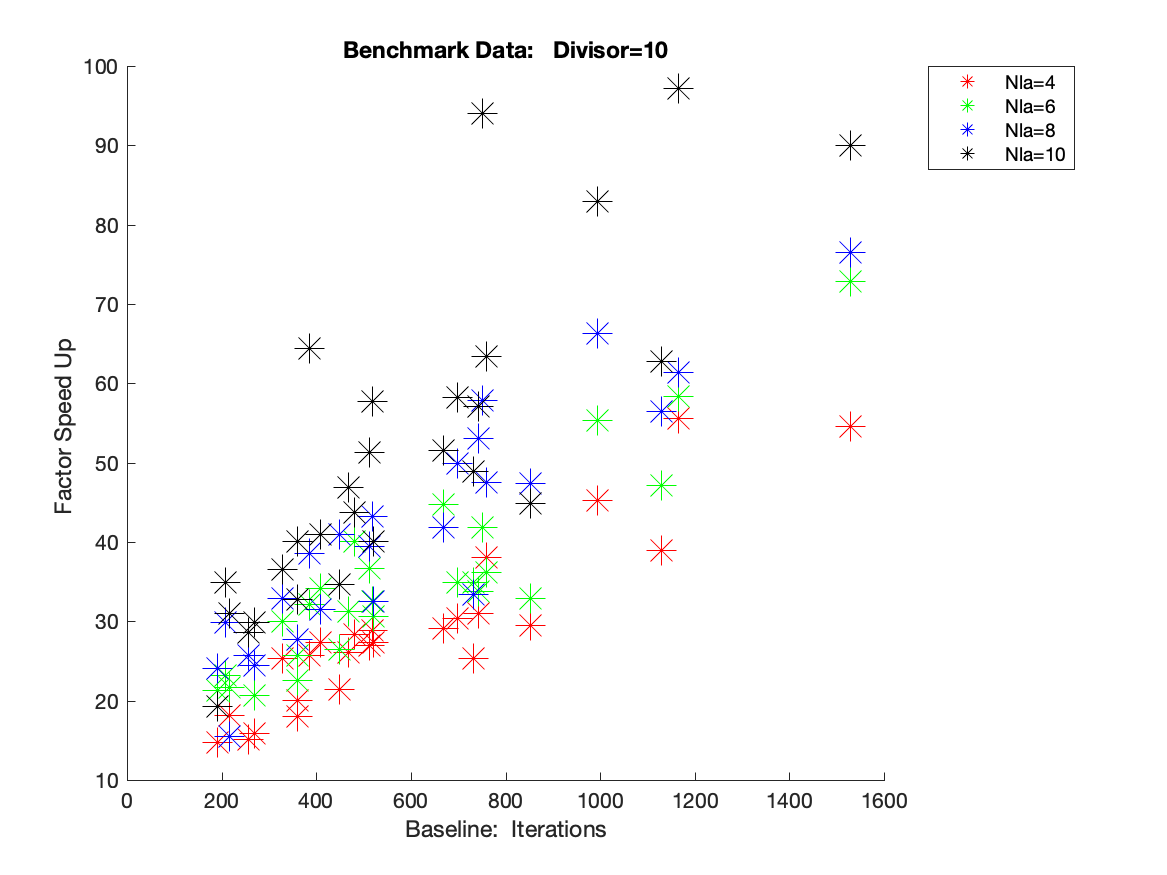}
\includegraphics[width=0.5\linewidth]{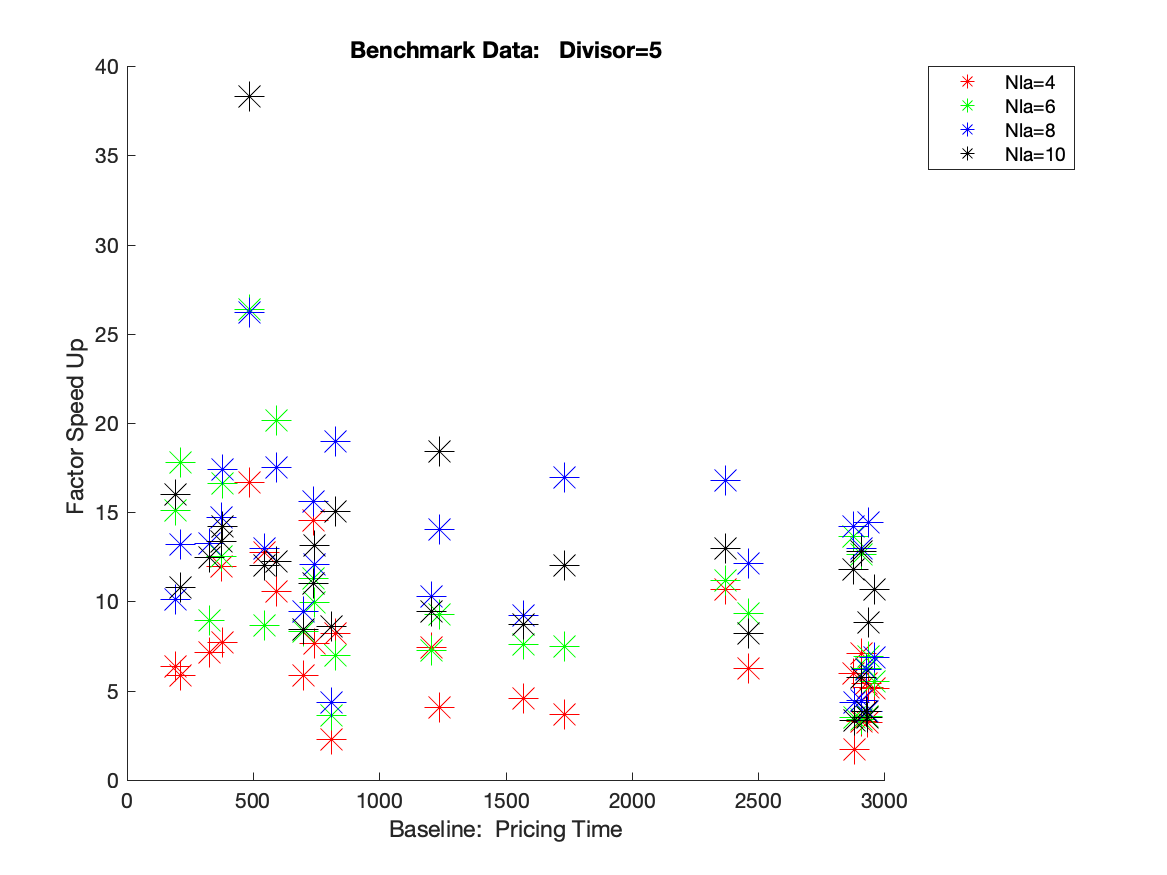}
\includegraphics[width=0.5\linewidth]{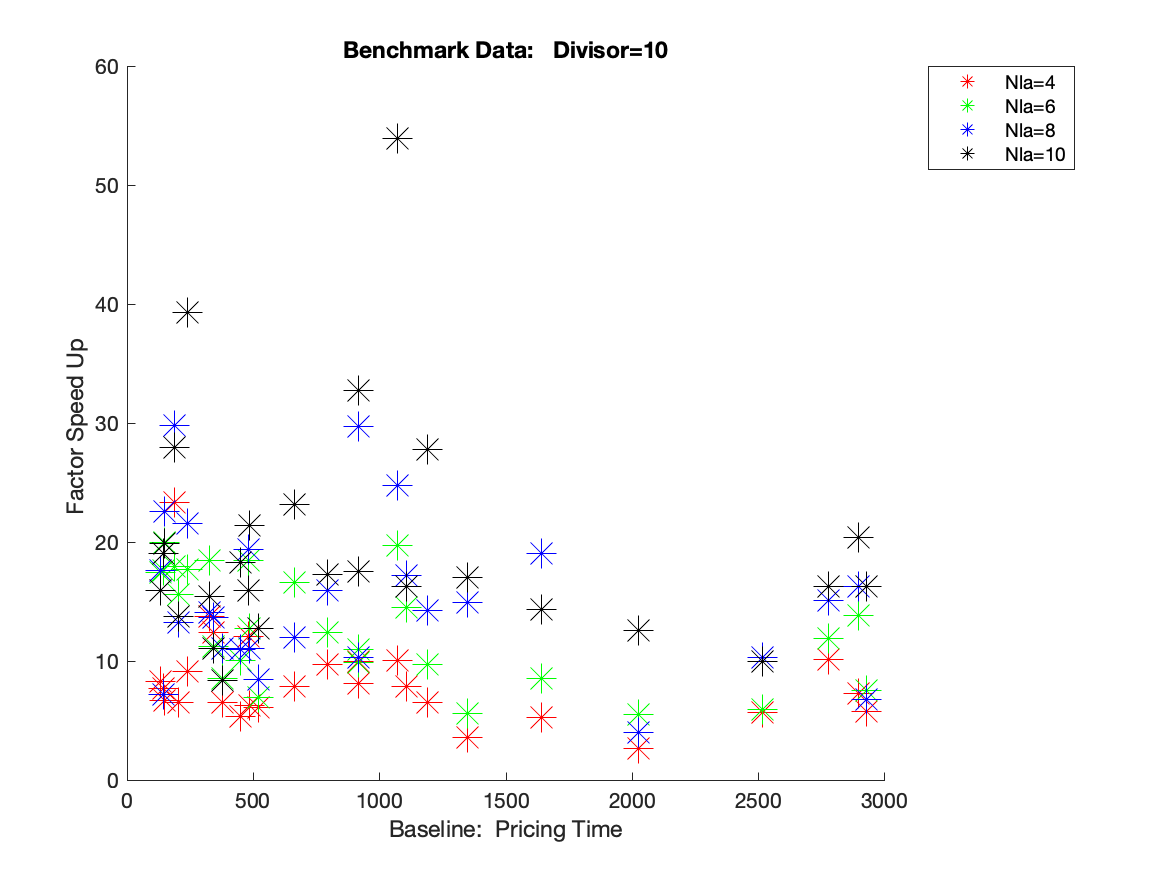}
\includegraphics[width=0.5\linewidth]{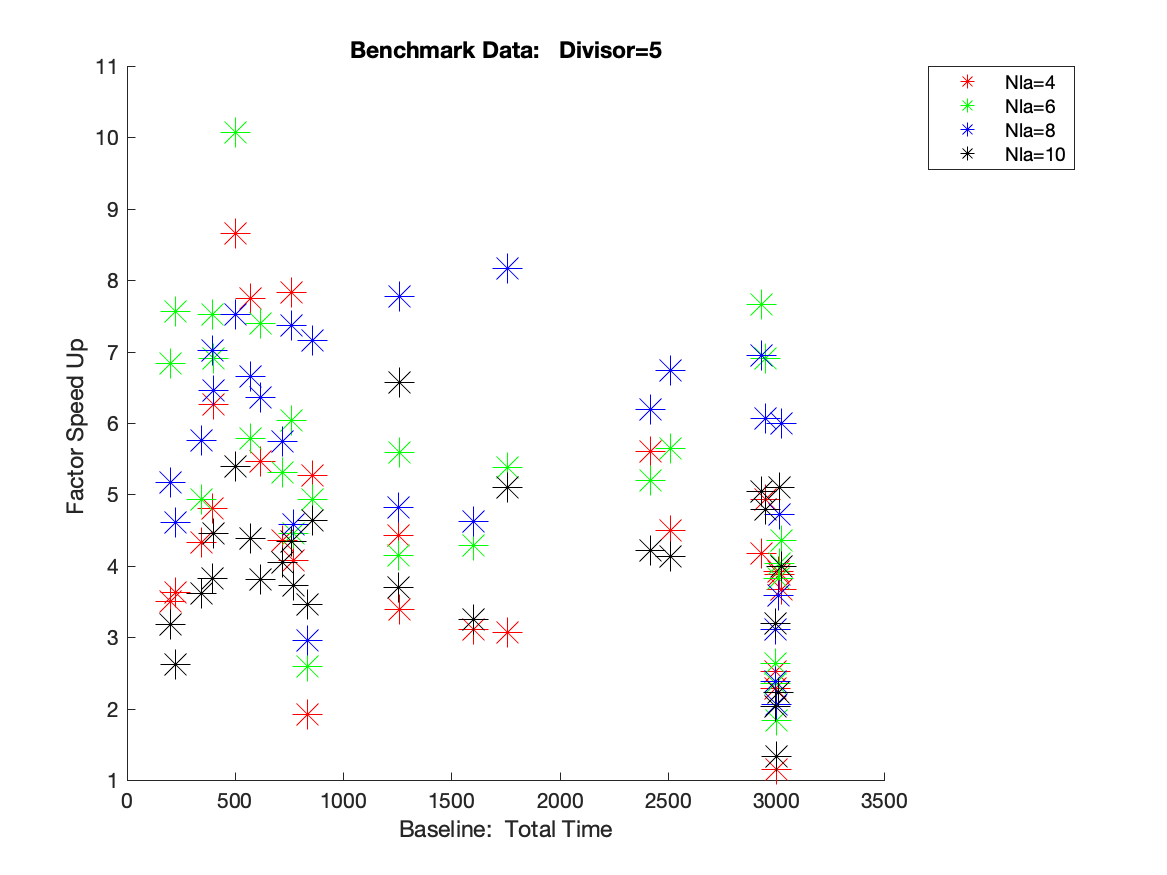}
\includegraphics[width=0.45\linewidth]{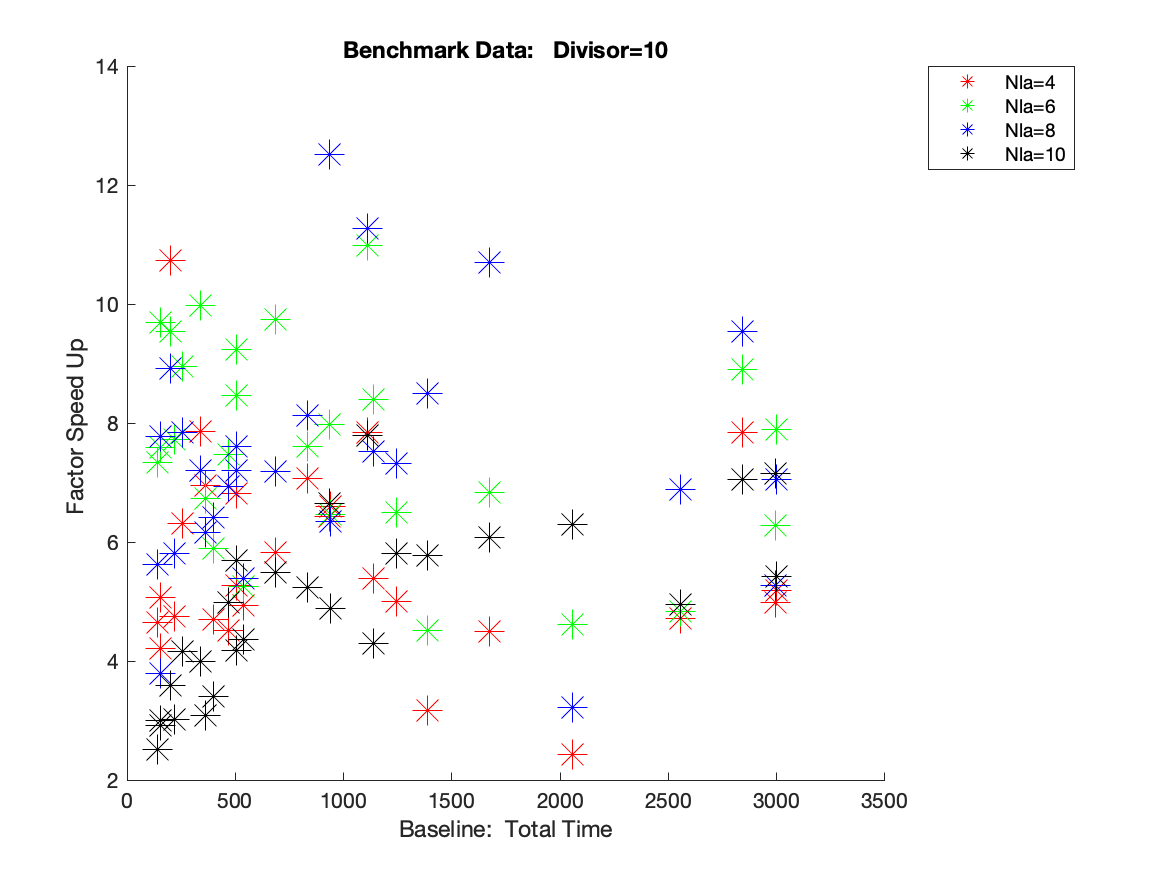}
\caption
{
\textbf{Top Row:} Relative performance of our CG Solver in terms of Iterations\\
\textbf{Middle Row:}  Relative performance of our CG Solver in terms of Pricing Time\\
\textbf{Bottom Row:} Relative performance of our CG Solver in terms of Total Time
    } 
\label{myCompTime_Augerat}
\end{minipage}
\end{figure*}
%
%\subsection{A Natural Question: What about Rounded %Capacity Inequalities?  Will LA-Routes Help there too?}
%Up to this point we have been concerned with implementing LA-Routes with Subset Row Inequalities (SRI), but Rounded Capacity Inequalities (RCI) are equally as valid and popular for CVRP problems. 
% See for example \citep{letchford2007branch, baldacci2012recent, diarrassouba2017complexity, diarrassouba2017separation}.  
% For completeness we examined the performance of LA-Route relaxation with RCI on the Augerat A problems. Their performance, along with a brief discussion of RCI is presented in the Appendix in Section \ref{RCI_Augerat}. 
%
%
%
\clearpage

\section{Conclusion}
\label{sec_conc}
We introduce Local Area (LA) route relaxations to improve the tractability and speed of Column Generation (CG) based solvers.  % for large scale set cover/partitioning formulations, where pricing is a elementary resource constrained shortest path problem \citep{costa2019}; this is a common framework found in large-scale transportation problems\citep{barnprice,Desrochers1992}.
We adapt subset row inequalities (SRI) \citep{jepsen2008subset} to efficiently tighten the ELP for VRP, producing LA-SRI. We integrate LA-SRI alongside LA-routes in a manner that ensures that the structure of the pricing problem is not altered. We demonstrate that our formulation with LA-SRI allows for an accelerated solution. We apply our approach to the CVRP, though our approach can be applied to other VRP \citep{Desrochers1992,costa2019}. We demonstrate that the ELP is indeed significantly tightened using LA-SRI and that computation of an optimal LP solution remains tractable.
%In future work, we intend to explore the use of time windows, other classes of valid inequalities in our formulations and exploit the primal block angular structure of $\Psi^+(\Omega_R,\Delta_R)$ for efficient solutions.
In future work we seek to use time windows in our formulation. This may involve different LA neighbor settings depending on the time when the vehicle leaves a customer. We also intend to use the original SRI in our formulation when they are violated by the LP solution but the associated LA-SRI are not. We shall seek to develop an algorithm that constructs the LA neighbor sets for customers in an efficient manner so as to maximally tighten the LP. In this process, we would aim to maximally tighten the LP bound given that the number of LA neighbors per customer does not exceed a user defined value (thus describing a computation budget). We shall further look into updating the heuristic used in A* search over the course of iterations of DSSR to avoid excessive expansion of nodes.  
 %Trivially we could compute for each possible tuple of [start and end time, LA-arc] the lowest cost time window feasible path. This path starts and ends at the same customers as the LA-arc and visits the associated customers s.t. the path starts at the start time and ends prior to the end time.  Then we could apply LA-routes on a graph where nodes describe the tuple of [remaining capacity, current time, ng-neighbors visited, current location]. More efficient mechanisms than enumerating all such time windows for each LA-arc would be a fruitful area of exploration.  Specifically 

%These adaptations will be done to allow the class of valid inequalities considered to easily fit into the existing pricing procedure, making the addition of such valid inequalities computationally efficient. 
%

\section*{Acknowledgement}

The authors thank Professor Louis-Martin Rousseau of the Polytechnique Montréal for his helpful suggestions. Any errors or omissions are those of the authors alone.

\bibliographystyle{abbrvnat} % outcomment this and next line in Case 1
%\singlespacing
\bibliography{col_gen_bib_2}
%
%\end{document}
% Acknowledgments here

% Appendix here
% Options are (1) APPENDIX (with or without general title) or 
%             (2) APPENDICES (if it has more than one unrelated sections)
% Outcomment the appropriate case if necessary
%
% \begin{APPENDIX}{<Title of the Appendix>}
% \end{APPENDIX}
%
%   or 
%
% \begin{APPENDICES}
% \section{<Title of Section A>}
% \section{<Title of Section B>}
% etc
% \end{APPENDICES}

% References here (outcomment the appropriate case) 

% CASE 1: BiBTeX used to constantly update the references 
%   (while the paper is being written).
%\bibliographystyle{informs2014trsc} % outcomment this and next line in Case 1
%\bibliography{<your bib file(s)>} % if more than one, comma separated

% CASE 2: BiBTeX used to generate mypaper.bbl (to be further fine tuned)
%\input{mypaper.bbl} % outcomment this line in Case 2
\newpage
\appendix
%\begin{APPENDIX}{}

In this appendix we discuss additional details about LA-Routes, LA-SRI and for completeness we provide an examination of the benefit of LA-SRI when combined with RCI. %In Appendix \ref{appendix_imp_details}, we discuss an implementation detail to improve the speed of our pricing procedure generating LA-routes. 
In Section \ref{sec_diff_LASRI_sri} we contrast LA-SRI and SRI, with regards to tightening the weighted set cover LP relaxation. In Section \ref{subseqEq} we demonstrate the equivalence of $\Psi^+(\Omega_R,\Delta_R)$ and $\Psi(\cup_{l \in \Omega_R}\Omega_l,\Delta_R)$.  In Section \ref{subsec_eff_solve}  we discuss the efficient solution of $\Psi^+(\Omega_R,\Delta_R)$.  %why we can use one equation in place of another to solve the Restricted Master Problem (RMP) for the Capacitated Vehicle Routing Problem (CVRP).
%In Section \ref{subsec_eff_solve} we discuss how we can provide an efficient solution to the RMP to CVRP while incorporating LA-routes and LA-SRI.  
In Section \ref{secAddExper}, we discuss additional timing results for problems with non-unit demand. In Section \ref{RCI_Augerat} we consider the use of rounded capacity inequalities (RCI) in the context of our solver.  % In Section \ref{other_ext} we discuss how we can improve upon our existing solver by improving LA neighbor set generation and other valid inequalities. % In Sections \ref{addLit1},\ref{addLit2} we consider additional literature review with regard to general dual stabilization and dual optimal inequalities respectively.
%
%\section{Implementation Details for LA-Route Pricing}
%\label{appendix_imp_details}
  %Thus this into pricing by not expanding such nodes $u,M_2,d_2$ when a corresponding node $u,M_1,d_1$ for $u,M_1\subseteq M_2\subseteq M_u,d_u\leq d_1\leq d_2$  where $u,M_1,d_1$ has already been expanded can improve solution time.
%
%
%
\section{Tightness of LA-SRI Relative to SRI}
\label{sec_diff_LASRI_sri}
In this section we compare tightening effect of a given SRI compared to the associated LA-SRI. % We show when an LA-SRI is weaker relative to the associated SRI. %Note that \ref{my_div_3}, is just a repeat of \ref{my_div_2} but needed for clarity here. 
%\begin{align}
%\label{my_div_3}
%\sum_{l\in \Omega}\theta_l \left\lfloor\frac{\sum_{u \in N_{\delta}}a_{ul}}{m_{\delta}}\right\rfloor\leq \left\lfloor \frac{|N_{\delta}|}{m_{\delta}}  \right\rfloor \; \forall N_{\delta}\subseteq N,  m_{\delta} \in \mathbb{Z}_{+}
%\end{align}
%
We observe with well chosen and large LA-neighborhoods that the relaxation of SRI to LA-SRI does not weakened the corresponding inequalities dramatically in practice. We provide substantial evidence for these phenomena with an analysis of the examples below; and define these examples over the following problem domain.

In our examples, suppose that we have a depot in San Diego (SD), and vehicles of capacity $201$.  We have three customers in New York City (NYC), named NYC1, NYC2, and NYC3, and thee customers in San Diego, named SD1, SD2, SD3. NYC and SD are far from each other but the associated customers of each city are nearly co-located. Each NYC customer has a demand of 100 and the SD customer has a demand of 1. Each customer has neighbor sizes of 2, so each customer in NYC considers one another to be an LA neighbor, and each customer in San Diego considers one another to be an LA neighbor. This is because we set LA neighbor sets by closest location.
%\begin{enumerate}
%
%\item 
%Consider a case with SRI defined to be $N_{\delta}$=[NYC1,NYC2,SD1] and $m_{\delta}$=2. Now consider routes $l_1$ and $l_2$. Route $l_1$ has the sequence of customers (with depot added ) as follows [-1, NYC1, NYC3, SD1, SD2, -2] and route $l_2$ has the sequence [-1, NYC1, NYC2, SD2, SD3, -2]. Note that the special indexes are 1,3 for both $l_1$ and $l_2$. Observe that for route $l_1$ the associated variables are defined to be $a_{\delta l_1}=1$ and $a^*_{\delta l_1}=0$. However, for route $l_2$ observe that $a_{\delta l_2}=1$ and $a^*_{\delta l_2}=1$. 
%\item 
%Consider a case with SRI defined to be $N_{\delta}$=[NYC1,NYC2,NYC3] and $m_{\delta}$=2. Now consider routes $l_3$ and $l_4$. Route $l_3$ has the sequence [-1, NYC1, NYC2, SD1, -2] and route $l_4$ has the sequence [-1, NYC1, SD1, NYC2, -2]. Note that the special indexes are 1,3 for $l_3$ and 1,2,3 for $l_4$. Observe that for route $l_3$ the associated variables are defined to be $a_{\delta l_3}=a^*_{\delta l_3}=1$. However, for route $l_4$ observe that $a_{\delta l_4}=1$ and $a^*_{\delta l_4}=0$. 
%
%\end{enumerate}
%
When solving the CVRP LP relaxation in \eqref{myRMPSSTAB} with the LA neighborhoods defined above alongside all LA-SRI of $|N_{\delta}|=3,m_{\delta}=2$, we enforce that at least two routes visiting at least one customer in NYC are used. %The two routes must visit different customers in SD and NYC.
When solving the CVRP LP formulation with LA neighborhoods and arcs but without LA-SRI, there is no such constraint on routes used. Given our relaxed formulation of LA-SRI constraints, one would think that larger LA neighborhoods always tightens the bound at the expense of greater computational cost during iterations of DSSR in pricing. However, we find that increasing the size of LA neighborhoods may not always tighten the bound, which we demonstrate with an example below.  

Consider that we add to each of the LA neighborhoods for (the three) customers in SD so that the LA neighbors of customers in SD are defined as follows:   $N_{SD1}=\{SD2,SD3,NYC1\}$,$N_{SD2}=\{SD1,SD3,NYC2\}$,$N_{SD3}=\{SD1,SD2,NYC3\}$. These LA neighborhood sets are still defined using closest distance. Now, consider the following solution to the MP $\Psi(\Omega,\Delta)$ using $\theta_{l_1}=\theta_{l_2}=\theta_{l_3}=\theta_{l_4}=\frac{1}{2}$ where routes $l_1,l_2,l_3,l_4$ are defined as sequences of the following LA-arcs.  

$l_1=\{(-1,SD1),(SD1,NYC1,NYC2),(NYC2,-2)\}$ \quad \\ $l_2=\{(-1,SD2),(SD2,NYC2,NYC3),(NYC3,-2)\}$,\quad 
\\$l_3=\{(-1,SD3),(SD3,NYC3,NYC1),(NYC1,-2)\}$,\quad 
\\$l_4=\{(-1,SD1),(SD1,SD2,SD3,-2)\}$

Observe that this solution uses $\frac{3}{2}$ cross country routes but no integer solution can use fewer than $2$ cross country routes.  

Based on the examples above we hypothesize that LA-SRI defined over customers localized in space are crucial for tightening the ELP relaxation of CVRP (and other VRP).  Furthermore, we observe that it is unlikely for an optimal solution to the  ELP relaxation to use routes that frequently return to areas of customers localized in space after leaving such areas. This is because including such routes in an LP solution would lead to sub-optimal cost. An example of such an route is [-1, NYC1, SD1, NYC2, -2]; as instead, the customers can be serviced at lower cost using the route [-1, NYC1, NYC2, SD1, -2].
Thus, for violated SRI defined over customers localized in space, a key opportunity for the corresponding LA-SRI to not be violated is when LA-arcs exist for which the following holds: the penultimate customer of the LA-arc is nearby the final customer of the LA-arc. By avoiding the construction of LA neighborhoods that allow for presence of such LA-arcs, we diminish the possibility of weakening the LP relaxation using LA-SRI instead of SRI. 
\section{Restricted Master Problem Equivalence}
\label{subseqEq}
In this section we establish that $\Psi^+(\Omega_R,\Delta_R)=\Psi(\cup_{l \in \Omega_R}\Omega_l,\Delta_R)$.  We rewrite $\Psi^+(\Omega_R,\Delta_R)$ for convenience below.  
\begin{subequations}
\label{myRMPSSTAB_A}
\begin{align}
    \min_{\substack{x\geq  0\\\theta \geq 0}}\sum_{l \in \Omega_R}c_l\theta_l +
    \sum_{\substack{l \in \Omega_R\\y \in Y^l\\ p \in \Omega^l_y}}c_{p}x^l_p+\sum_{\substack{l \in \Omega_R\\ u \in N\\ d_0\geq d\geq d_u}}c_{-1 u}x^l_{(-1,d_0),(u,d)}
    \label{RMP_stab_obj_A}\\
    \sum_{l \in \Omega_R}\theta_l+
    \sum_{l \in \Omega_R}\sum_{ij \in E^l}[i=(-1,d_0)]x^l_{ij}\leq K\quad [-\pi_0] \label{RMP_stab_packing_A}\\
\sum_{l \in \Omega_R}a_{ul}\theta_l+
\sum_{\substack{l \in \Omega_R\\y \in Y^l\\ p \in \Omega^l_y}}x^l_{p}a_{up}\geq 1 \quad \forall u \in N\quad [\pi_u] \label{RMP_stab_cover_A}\\
    \sum_{l \in \Omega_R}a^*_{\delta l}\theta_l+
    \sum_{\substack{l \in \Omega_R\\y \in Y^l\\ p \in \Omega^l_y}}a^*_{\delta p}x^l_p\leq \left \lfloor \frac{|N_{\delta}|}{m_\delta}\right \rfloor \quad \forall \delta \in \Delta_R \quad [-\pi_{\delta}]\label{RMP_stab_delta_A}\\
    \sum_{p \in \Omega^l_y}x^l_{p}= \sum_{\substack{ij\in E^l\\i=(u,d_1)\\ j=(v,d_1-d)}}x^l_{ij} \label{RMP_stab_agree_A} \forall y \in Y^l; y=(u,v,d),l \in \Omega_R\\
    \sum_{ij \in E^l}x^l_{ij}=\sum_{ji \in E^l}x^l_{ji} \quad \forall i \in I^0,  u_i\notin \{-1,-2\},l \in \Omega_R \label{RMP_stab_flow_A}
\end{align}
\end{subequations}

We establish that in two parts. In Section \ref{subsec_first} we demonstrate that $\Psi^+(\Omega_R,\Delta_R)\leq \Psi(\cup_{l \in \Omega_R}\Omega_l,\Delta_R)$.  In Section \ref{subsec_second} we establish that $\Psi^+(\Omega_R,\Delta_R)\geq \Psi(\cup_{l \in \Omega_R}\Omega_l,\Delta_R)$.  Since $\Psi^+(\Omega_R,\Delta_R)\leq \Psi(\cup_{l \in \Omega_R}\Omega_l,\Delta_R)$ and $\Psi^+(\Omega_R,\Delta_R)\geq \Psi(\cup_{l \in \Omega_R}\Omega_l,\Delta_R)$  then $\Psi^+(\Omega_R,\Delta_R)= \Psi(\cup_{l \in \Omega_R}\Omega_l,\Delta_R)$.
\subsection{First Side of the Inequality}
\label{subsec_first}
In this section we  establish that $\Psi^+(\Omega_R,\Delta_R)\geq \Psi(\cup_{l \in \Omega_R}\Omega_l,\Delta_R)$. % using proof by contradiction.   
Consider the optimal solution to \eqref{myRMPSSTAB_A} denoted $(x,\theta)$. 
Observe that if the $x$ terms are zero valued then $\Psi^+(\Omega_R,\Delta_R)=\Psi(\Omega_R,\Delta_R)$ and $\Psi(\Omega_R,\Delta_R)\geq  \Psi(\cup_{l \in \Omega_R}\Omega_l,\Delta_R)$ by definition thus establishing the claim.  

We now transform the solution $x,\theta$ to a solution to  $\Psi(\cup_{l \in \Omega_R}\Omega_l,\Delta_R)$ with identical cost to that of $x,\theta$. Each step of this process decreases the number of $x$ terms with non-zero value.  We now describe an individual step of this process.  If $x$ is not zero valued by \eqref{RMP_stab_flow_A}, there must exist a path of non-zero valued $x^l_{ij}$ terms starting at the source node $(-1,d_0)$ and ending at the sink node $(-2,0)$ for some $l \in \Omega_R$.  
Select any such path on the graph associated with $l$, and let the edges on the path be denoted $\hat{E}^l$.  By \eqref{RMP_stab_agree_A} for each $ij\in \hat{E}^l$ (except the edge including the source) there must exist a $p \in \Omega^l_y$ (for $i=(u,d_1)$,$j=(v,d_2)$, $y=(u,v,d_1-d_2$) s.t. $x^l_p>0$.  For each $ij \in \hat{E}^l$ (excluding the edge connected to the $(-1,d_0)$) let $p_{ij}$ be any such $p \in \Omega^l_y$ s.t. $x^l_p>0$. 
Let $\hat{l}$ denote a route corresponding to crossing all $ij \in \hat{E}^l$ and using the paths of the associated $p_{ij}$ as the LA-arcs describing intermediate customers visited.  Note that $\hat{l}$ is an elementary feasible route that lies in $\Omega_l$ (but need not lie in $\Omega_R$) by definition of $E^l,Y^l$.  Now consider the following altered solution given a tiny constant $\alpha>0$.  
\begin{subequations}
\label{my_grad_step}
\begin{align}
    \theta_{\hat{l}} \leftarrow \theta_{\hat{l}} +\alpha\\
    x^l_{ij} \leftarrow x^l_{ij} -\alpha \quad \forall ij \in \hat{E}^l\\
    x^l_{p}\leftarrow x^l_p- \alpha \quad \forall p \in \cup_{\substack{ij \in \hat{E}^l\\ i\neq (-1,d_0)}}p_{ij}
\end{align}
\end{subequations}
Observe that the change in cost is zero by the definition of the cost of a route.
Thus we need to establish that a non-zero value $\alpha$ exists s.t. the change induced in \eqref{my_grad_step} is feasible.  Below we set $\alpha$ to the largest possible value s.t. \eqref{my_grad_step} does not alter non-negativity.% which is defined as follows:
\begin{align}
\label{alphaDed}
    \alpha \leftarrow  \min \{ \min_{ij \in \hat{E}^l}x^l_{ij},\min_{\substack{ij \in \hat{E}^l\\ i\neq (-1,d_0)}}x^l_{p_{ij}})\}
\end{align}
Since all $x$ terms in \eqref{alphaDed} are positive then $\alpha$ is positive. Observe that this solution has fewer non-zero valued $x$ terms than prior to the update. 
\subsection{Second Side of the Inequality}
\label{subsec_second}
In this section we establish that $\Psi^+(\Omega_R,\Delta_R)\leq \Psi(\cup_{l \in \Omega_R}\Omega_l,\Delta_R)$ using proof by contradiction.  Consider any solution producing optimal terms $\theta$ to $ \Psi(\cup_{l \in \Omega_R}\Omega_l,\Delta_R)$. %
%If $\Psi^+(\Omega_R,\Delta_R)\geq \Psi(\cup_{l \in \Omega_R}\Omega_l,\Delta_R)$ then a contradiction is created establishing the claim.  %established.  

We construct as solution $x$ (with zero valued $\theta$) as follows to $\Psi^+(\Omega_R,\Delta_R)$ with the same objective as that of $\theta$ over $\Psi(\cup_{l \in \Omega_R}\Omega_l,\Delta_R)$.  Let $\Omega_{R2}=\cup_{l \in \Omega_R}\Omega_l$.  For for each $l\in \Omega_{R2}$ let $Q_l$ be any route in $\Omega_R$ for which  $l \in \Omega_{Q_l}$.  %Now iterate of $l\in \Omega$ for which $\theta_l>0$%(where $\Omega_{Q_l}$ is simply $\Omega_{\hat{l}}$ where $\hat{l}=Q_l$).
We now construct $x$ as follows.  
\begin{subequations}
\begin{align}
   x_{ij}^{\hat{l}}\leftarrow \sum_{\substack{l\in \Omega_{R2}\\ \hat{l}=Q_l}} a_{ijl}\theta_l \quad \forall (ij \in E^{\hat{l}}, \hat{l} \in \Omega_R)\\
   x_{p}^{\hat{l}}\leftarrow \sum_{\substack{l\in \Omega_{R2} \\ \hat{l}=Q_l}} a_{pl}\theta_l \quad \forall p \in \Omega^{\hat{l}}_y, y\in Y^{\hat{l}},\hat{l} \in \Omega_R
\end{align}
\end{subequations}
Observe that this solution $x$ (with $\theta$ set to zero) is feasible and has cost identical to $\Psi(\cup_{l \in \Omega_R}\Omega_l,\Delta_R)$ creating a contradiction hence establishing the claim.  
%
%\end{document}

\section{Efficient Solution to the Restricted Master Problem}
\label{subsec_eff_solve}
In this section we consider the efficient solution of \eqref{myRMPSSTAB}.  Observe that \eqref{myRMPSSTAB} may have a very large number of variables as $|\Omega_R|$ grows and for problems where $|\Omega_y|$ is large for some $y \in Y$; thus making the solution to \eqref{myRMPSSTAB} at each iteration of CG intractable.  %Hence we solve \eqref{myRMPSSTAB} in a manner akin to CG.
Thus we seek to construct a small set of edges, paths  denoted $\hat{E}^l\subseteq E^l \quad \forall l \in \Omega_R,\quad \hat{\Omega}_{y}^l \subseteq \Omega_{y}^l$ \quad $(y \in Y^l \; \forall l \in \Omega_R)$ respectively s.t. solving \eqref{myRMPSSTAB} over these terms (denoted $\Psi^+(\Omega_R,\Delta_R,\{ \hat{E}^l\},\{\hat{\Omega}^l_y\}))$ yields the same solution as \eqref{myRMPSSTAB};  where for short hand we use $\{ \hat{E}^l \}, \{ \hat{\Omega}^l\} $ to describe the $\hat{E}^l,\hat{\Omega}^l_y$ for each $l \in \Omega_R,(y\in Y^l,l \in \Omega_R)$.  
We define $\Psi^+(\Omega_R,\Delta_R, \{ \hat{E}^l\} ,\{ \hat{\Omega}^l_y \} ))$ formally below.

\begin{subequations}
\label{subPartEqRMP}
\begin{align}
    %\Psi^+(\Omega_R,\Delta_R,\{ \hat{E}^l \} ,\{ \hat{\Omega}^l_y\} )=\\
    \min_{\substack{x\geq  0}}
    \sum_{\substack{l \in \Omega_R\\y \in Y^l\\ p \in \hat{\Omega}^l_y}}c_{p}x^l_p+\sum_{\substack{l \in \Omega_R\\ u \in N\\ d_0\geq d\geq d_u\\ (-1,d_0),(u,d) \in \hat{E}^l}}c_{-1 u}x^l_{(-1,d_0),(u,d)}\\
    \sum_{l \in \Omega_R}\sum_{ij \in \hat{E}^l}[i=(-1,d_0)]x^l_{ij}\leq K\quad [-\pi_0] \\
\sum_{\substack{l \in \Omega_R\\y \in Y^l\\ p \in \hat{\Omega}^l_y}}x^l_{p}a_{up}\geq 1 \quad \forall u \in N\quad [\pi_u] \\
    \sum_{\substack{l \in \Omega_R\\y \in Y^l\\ p \in \hat{\Omega}^l_y}}a^*_{\delta p}x^l_p\leq \frac{|N_{\delta}|}{m_{\delta}} \quad \forall \delta \in \Delta_R \quad  [-\pi_{\delta}]\\
    \sum_{p \in \hat{\Omega}^l_y}x^l_{p}= \sum_{\substack{ij\in \hat{E}^l\\i=u,d_1\\ j=v,d_1-d}}x^l_{ij}\quad \forall y \in Y^l; y=(u,v,d),l \in \Omega_R \\
    \sum_{ij \in \hat{E}^l}x^l_{ij}=\sum_{ji \in \hat{E}^l}x^l_{ji} \quad \forall i\in I^0, \quad u_i\notin \{-1,-2\},l \in \Omega_R \quad 
\end{align}
\end{subequations}
To solve \eqref{myRMPSSTAB} we alternate between the following two steps.  
\begin{itemize}
    \item Solve $\Psi^+(\Omega_R,\Delta_R,\{ \hat{E}^l \}, \{ \hat{\Omega}^l_y \} )$ producing the dual solution  $\pi$.  This is fast since the sets $\hat{\Omega}^l_y$ and $\hat{E}^l$ are generally much smaller than $\Omega^l_y, E^l$ respectively. % The output dual variables are denoted $\pi$.  
    \item  Iterate over $l \in \Omega_R$ and compute the lowest reduced cost route in $\Omega_l$ denoted $\hat{l}^l$.  The computation of $\hat{l}^l$ is a simple shortest path computation (not a resource constrained shortest path problem) and corresponds to the lowest cost path on the graph with edge set $E^l$ where the edge weights are defined as follows. 
    \begin{subequations}
    \label{my_weight_eq}
    \begin{align}
    \bar{c}_{(-1,d_0),(u,d)}=c_{-1,u}+\pi_0\\
    \bar{c}_{ij}=\min_{\substack{p \in \Omega^l_y\\ y=(u,v,d_1-d_2)}}\bar{c}_p \\\quad \forall i=(u,d_1),j=(v,d_2), i,j \in E^l \nonumber \\
    \bar{c}_p=c_p-\sum_{u \in N}a_{up}\pi_u+\sum_{\delta \in \Delta_R}a^*_{\delta p}\pi_{\delta} \\ \forall p \in \Omega^l_y, y\in Y^l\nonumber 
    \end{align}
    \end{subequations}
    Observe that via the ordering of $\beta^l$ that the route $\hat{l}^l$ must be elementary as all paths from source to sink describe elementary routes.  
    %We then add to $\hat{E}^l,\hat{\Omega}^l_y$ (for the relevant $y \in Y^l$) all edges, paths used in $\hat{l}^l$.
    If $\hat{l}^l$ has negative reduced cost then we update the edges,paths associated with route $l$.  We denote the edges used in this route as $E^l(\hat{l})$ and the paths used in this route for a given $y$ in $Y^l$ as $\Omega^l_y(\hat{l})$.  Next we augment each $\hat{E}^l,\hat{\Omega}^l_y (\forall y \in Y^l)$ with the edges with $E^l(\hat{l}^l)$ and $\Omega^l_y(\hat{l}^l)$.   When we find that $\hat{l}^l$ has non-negative reduced cost for each $l\in \Omega_R$, we terminate since we have solved \eqref{myRMPSSTAB} optimally.
    
    We refer to the operation of generating edge weights and computing the associated shortest path as the RMP-Shortest Path computation (RMP-SP).  
\end{itemize}
We must initialize the $\{ \hat{\Omega}^l_y \},\{ \hat{E}^l\} $  terms in order to ensure that \eqref{subPartEqRMP} has a feasible solution. In experiments we found that using $\hat{\Omega}^l_y,\hat{E}^l$ to be the terms generated during RMP-SP for a route $l$ over the entire course of CG optimization thus far worked well. We could just as easily use all edges that are associated with non-zero $x^l_{ij},x^l_{p}$ values in the final solution produced last time \eqref{myRMPSSTAB} was solved using \eqref{subPartEqRMP}.
In our experiments (which had no maximum number of vehicles) prior to the first iteration of CG we included all edges and paths corresponding to the route starting at the depot, visiting a single customer $u$, and then returning to the depot (for each $u\in N$). 
This means that for the initial $l$ in CG, denoted $l_0$ we define $\hat{E}^{l_0}$ to include $(-1,d_0),(u,d_u)$,and $(u,d_u),(-2,0)$ for each $u \in N$; and set $\hat{\Omega}^{l_0}_y$ as empty except for terms $y=(u,-2,d_u)$ (for each $u\in N$) which is set to contain the path starting at $u$ ending at $-2$ and visiting no intermediate customers.  

In Alg \ref{fancyPGM} we provide a formal description of our fast optimization of \eqref{myRMPSSTAB} using RMP-SP which we annotate below.  %and we provide annotations below.  

\begin{algorithm}[!b]%[1]%[H]
 \caption{Fast Solver for \eqref{myRMPSSTAB}}
\begin{algorithmic}[1] 
\State $\Omega_R\{ \hat{E}^l\} ,\{ \hat{\Omega}_y^l\} \leftarrow $ from user \label{init_1}
\Repeat \label{start_loop}%{True} 
\State $x,\pi \leftarrow  $ Solve \eqref{subPartEqRMP} over $\hat{E}^l,\hat{\Omega}_y^l$ \label{solve_PGM_RMP}
\For{$l \in \Omega_R$}\label{start_l}
\State $\hat{l}^l\leftarrow  \mbox{arg} \min_{l \in \Omega_l}\bar{c}_l$.  Fast and easy shortest path calculation (not resource constrained).%Solve for lowest reduced cost $\hat{l}\in \Omega_l$ (not $\Omega$).   % This computation is a simple shortest path computation on a small graph with weights \eqref{my_weight_eq} %and one node for the su
\label{start_price_in_family}
\If{$\bar{c}_{\hat{l}^l}<0$} \label{add_family_check}
\State $\hat{E}^l \leftarrow \hat{E}^l \cup E^l(\hat{l}^l)$% \cup_{[ij]\in E^l}ij $ s.t. $ij$ is used in route $\hat{l}^l$ meaning $a_{ij\hat{l}^l}=1$ 
\label{add_in_family}
\State $\hat{\Omega}^l_y \leftarrow \hat{\Omega}^l_y  \cup \Omega^l_y(\hat{l}^l)$ for each $y \in Y^l$ 
\label{add_in_family2} 
\EndIf% \cup_{[ij]\in E^l}ij $ s.t. $ij$ is used in route $\hat{l}^l$ meaning $a_{ij\hat{l}^l}=1$ \label{add_in_family}
% \hat{\Omega}^l_y \cup_{\substack{p \in \Omega^l_y\\ a_{p\hat{l}^l}=1}}p $.
%,p \in \Omega^l_y$. % s.t. arc $p$ is used in route $\hat{l}_l$ meaning $a_{p\hat{l}_l}=1$.  Apply for all $p$ 
\EndFor \label{end_l}
 \Until{$\bar{c}_{\hat{l}^l}\geq 0$ for all $l \in \Omega_R$} \label{term_me} %a_{\delta_*l}\theta_l\geq 0$}
 \For{$l \in \Omega_R$}\label{start_add}
 \State $\hat{E}^l \leftarrow \{ij \in \hat{E}^l \quad s.t. \quad x^l_{ij}>0\}$ OPTIONAL NOT USED IN EXPERIMENTS
  \State $\hat{\Omega}^l_y \leftarrow \{ p \in \hat{\Omega}^l_y  \quad s.t. \quad x^l_{y}>0\}$ for each $y\in Y^l$.   OPTIONAL NOT USED IN EXPERIMENTS
 \EndFor \label{end_add}
 \State Return last $x,\pi ,\{ \hat{E}^l\}, \{ \hat{\Omega}^l_y \} $% generated and $ \{ \hat{E}^l \quad \forall l \in \Omega_R\}, \{ \hat{\Omega}^l_y\quad \forall l \in \Omega_R,y \in Y^l\} $    %and $\hat{Q}^l$.% The terms for which $\psi^l_{uvd}>0$ become the new inital $\hat{Q}_l$ next time I solve RMP \label{returner}%Can be used inside branch-cut-price \label{ZreturnSol}
\end{algorithmic}
\label{fancyPGM}
\end{algorithm} 

\begin{itemize}
    \item In Line \ref{init_1} we initialize $\Omega_R,\{ \hat{E}^l\} ,\{ \hat{\Omega}^l_y\} $ from the user. 
    %There is one $\hat{E}^l$ for each $l \in \Omega_R$ and one $\hat{\Omega}^l_y$ for each $l \in \Omega_R,y \in Y^l$. %$E^l$ is the set of all edges associated with family $l$, $\hat{E}$
    \item In the loop defined in Lines \ref{start_loop}- \ref{end_add} we construct a sufficient set of $\{ \hat{E}^l\} ,\{ \hat{\Omega}^l_y\}$ to solve $\Psi^+(\Omega_R,\Delta_R)$ exactly meaning $\Psi^+(\Omega_R,\Delta_R,\hat{E},\hat{\Omega})=\Psi^+(\Omega_R,\Delta_R)$.
    \begin{itemize}
        \item In Line \ref{solve_PGM_RMP} we solve $\Psi^+(\Omega_R,\Delta_R,\{ \hat{E}^l\},\{ \hat{\Omega}^l_y \})$ as an LP producing $\pi,x$ terms. % This can be understood as solving 
        \item In Lines \ref{start_l}-\ref{end_l} we iterate over $l\in \Omega_R$ and augment  $\{ \hat{E}^l\},\{ \hat{\Omega}^l_y\}$ using RMP-SP.
        \begin{itemize}
            \item In Lines \ref{start_price_in_family} we generate the lowest reduced cost column in the family $\Omega_l$ which we denote as $\hat{l}^l$. It is important to note that this is a simple shortest path and not a hard resource constrained shortest path problem.
            \item In Line \ref{add_family_check} we determine if $\hat{E}^l,\hat{\Omega}^l_y$ need to be augmented which is true if $\hat{l}^l$ has negative reduced cost.  In Line \ref{add_in_family} we add any edges used in $\hat{l}^l$ to $\hat{E}^l$ expanding the set of edges.  In Line \ref{add_in_family2} we similarly, add any paths used in $\hat{l}^l$ to the corresponding $\hat{\Omega}^l_y$.%sets for each $y \in Y^l$ for $l$.
        \end{itemize}
        \item In Line \ref{term_me} we terminate the loop once no path in any graph $E^l$ for $l\in \Omega_R$ has negative reduced cost, since at this point we have solved $\Psi^+(\Omega_R,\Delta_R)$ exactly.
          \end{itemize}
  \item In Lines \ref{start_add}- \ref{end_add} we describe an optional step to ensure that $\hat{E}^l,\hat{\Omega}^l_y$ do not grow too large causing computational difficulty. In practice we found this unnecessary. In this step we remove all terms for which the corresponding  $x^l_{ij}$,$x^l_{p}$ are zero valued for the solution of \eqref{myRMPSSTAB} for the final solution to \eqref{subPartEqRMP}.    
\end{itemize}
% We descirbe teh soltuion 
%We initalize $\hat{\Omega}^l_y,\hat{E}^l$ with the elements used  In experiments in order to solve the RMP we seek to not use excessive numbers of calls.  Thus we never remove elements from  $\hat{E}^l,\hat{\Omega}{yl} \quad \forall l \in \Omega_R$.  Thus these are the set of edges and paths generated in any point in solving the RMP during shortest path calculation.  W
%The complete solution to \eqref{myRMPSSTAB} is thus written in Alg \ref{basicCGRLA} using the fast solution for the RMP \eqref{myRMPSSTAB} written in Alg \ref{fancyPGM}.  We provide annotation for these algorithms below.
Given Alg \ref{fancyPGM} we describe the complete optimization approach for solving $\Psi(\Omega,\Delta)$ using column/row generation in Alg \ref{basicCGRLA} with annotation below.  

\begin{algorithm}[!b]%[1]%[H]
 \caption{Complete Solver for the Master Problem}
\begin{algorithmic}[1] 
\State $\Omega_R,\leftarrow $ from user \label{Zline_rec_input_start1}
\State $\hat{E}^l,\hat{\Omega}_{y}^l$ from user for all $l \in \Omega_R$,$(y \in Y^l,l\in \Omega_R)$ respectively \label{Zline_rec_input_start0}
\item $\Delta_R \leftarrow \{ \}$
\label{delta_init}
\Repeat%{True}
\label{Zline_outer_start2a}%
\Repeat%{True}
\label{Zline_outer_start2}%
\State  $x,\pi,\{ \hat{E}^l \},\{ \hat{\Omega}^l_y\} \leftarrow $ Solve for \eqref{subPartEqRMP}, using Alg \ref{fancyPGM} \label{callBasicPGM} %\eqref{primal_master}
\label{Zline_solve_FRMP2}
\State $l_* \leftarrow \mbox{arg}\min_{l \in \Omega}\bar{c}_l$ via the LA-routes pricing algorithm \label{gen_col}
\State $\Omega_R \leftarrow \Omega_R \cup l_*$ \label{add_col}
%\State $\hat{E}_{l_*} \leftarrow \{ij \in E^{l_*}:  \quad s.t. \quad  a_{ijl_*}=1 \} $\label{definiteIntint1}
%\State $\hat{\Omega}^{l_*}_z \leftarrow \{ij \in \Omega^{l_*}_z:  \quad s.t. \quad  a_{pl_*}=1  \quad \}  \forall z \in Z^l$\label{definiteIntint2}% or $\hat{Q}^l=\{$
 \Until{$\bar{c}_{l_*} \geq 0$}\label{returnCond22}
 \State ${\delta}_*\leftarrow $ Most violated constraint over \eqref{RMP_stab_delta} given $x$. More than one violated constraint can be added.   \label{my_sep2}
 \State $\Delta_R \leftarrow \Delta_R \cup \delta_*$ \label{add_sri_line}
 \Until{\eqref{RMP_stab_delta} is not violated for $\delta_*$}\label{returnCond_last}
 \State Return last $x$  generated by line \ref{callBasicPGM}.  This $x$ can be used inside branch \& price\label{ZreturnSol2}.% to enforce integral solutions to $x$.
\end{algorithmic}
\label{basicCGRLA}
\end{algorithm}

\begin{itemize}
  %  \item IN PROGRESS
    \item In Lines \ref{Zline_rec_input_start1}-\ref{Zline_rec_input_start0} we initialize from the user the set of columns in $\Omega_R$. In our experiments this is a single column $l$ randomly generated with edge,path sets $\hat{E}^l,\hat{\Omega}^l_y \quad \forall y \in Y^l$ that can be used to generate all routes containing a single customer (one for each $u\in N$).% for all customers $u \in N$.
    %\item Line \ref{Zline_rec_input_start0}: Include all edges so that each route consisting of a single customer can be generated.  We set $\hat{E}^l,\hat{\Omega}^l_y$.
    \item In Line \ref{delta_init} we initialize the set $\Delta_R$ to be empty.
    %\item Line \ref{Zline_rec_input_start2} Define $\hat{E}_l$ to be the set of all $ij$ edges in family/column $l$ in $\Omega_R$
    %\item Line \ref{Zline_rec_input_start3} Define $\hat{\Omega}^l_z$ to be the set of all arcs ($u,v,M1,M2,d$) used in family/column $l$ in $\Omega_R$
    \item In Lines \ref{Zline_outer_start2a}-\ref{returnCond_last} we solve the master problem over all $\Omega,\Delta$, which is written $\Psi(\Omega,\Delta)$.
    \begin{itemize}
    \item In Lines \ref{Zline_outer_start2}-\ref{returnCond22} we solve $\Psi(\Omega,\Delta_R)$. %master problem given fixed set $\Delta_R$.
    \begin{itemize}
        \item In Line \ref{callBasicPGM} we solve $\Psi^+(\Omega_R,\Delta_R)$ using Alg \ref{fancyPGM}.  
        \item In Line \ref{gen_col} we find the lowest reduced cost route $l_*$.% Recall that  $\bar{c}_l$ is defined in \eqref{red_eq_1} for each $l \in \Omega$.
        %  Unlike in PGM (Alg \ref{fancyPGM}) this is a resource constrained shortest path problem.  We solve this by creating all possible LA-arcs using dynamic programming and finding the shortest path over these arcs adjusted using dual variables (note the extra $\pi_{\delta}$, the addition of such terms does not complicate the program).
        \item In Line \ref{add_col} we add the column $l_*$ to the $\Omega_R$ set.  
        %\item Line \ref{definiteIntint1}:  Define the initial edges in the $\hat{E}^{l_*}$ for the new route that I just added to $\Omega_R$
        %\item Line \ref{definiteIntint2}:  Define the initial $p$ in the $\hat{\Omega}^{l_*}_z$ for the new route that I just added to $\Omega_R$ .  This is done for each $z \in Z^{l_*}$
        \item In Line \ref{returnCond22} we terminate CG when the lowest reduced cost route in $\Omega$ has non-negative reduced cost. % This means that given fixed $\Delta_R$ that we have solved the master problem over $\Omega$ and $\Delta_R$ (NOT $\Delta$)
    \end{itemize}
    \item In Lines \ref{my_sep2}-\ref{add_sri_line} we find one or more violated LA-SRI which are then added to $\Delta_R$. In our experiments we add the 30 most violated LA-SRI (or all violated LA-SRI if less than 31 are violated)
    %\item  Line \ref{add_sri_line}: Add $\delta_*$ to $\Delta_R$.   % SRI set to th.e set of all SRI's we constrain over by using
    \item In Line \ref{returnCond_last} we terminate optimization when no LA-SRI are violated since we have solved $\Psi(\Omega,\Delta)$ optimally.
        \end{itemize}
    \item Line \ref{ZreturnSol2}:  Return $x$ providing a primal solution that can be fed into branch \& price.  
\end{itemize} 
It is often useful to generate an approximate optimal solution to the optimal integer linear programming (ILP) solution to CVRP.  To do this efficiently we  solve the RMP in $\Psi^+(\Omega_R,\Delta_R,\{ \hat{E}^l\},\{\hat{\Omega}^l\})$ as described in \eqref{subPartEqRMP}  as an ILP  given $\{ \hat{E}^l \}, \{ \hat{\Omega}^{l}_{y}\}$ generated during Alg \ref{basicCGRLA}. We found that this produces quality integer solutions in practice.
%, and these solutions can be used as an upper bound to the corresponding LP solution produced in \eqref{subPartEqRMP}. 
%

\section{Additional Experiments}
\label{secAddExper}

\begin{figure}[!hbtp]

\includegraphics[width=0.45\linewidth]{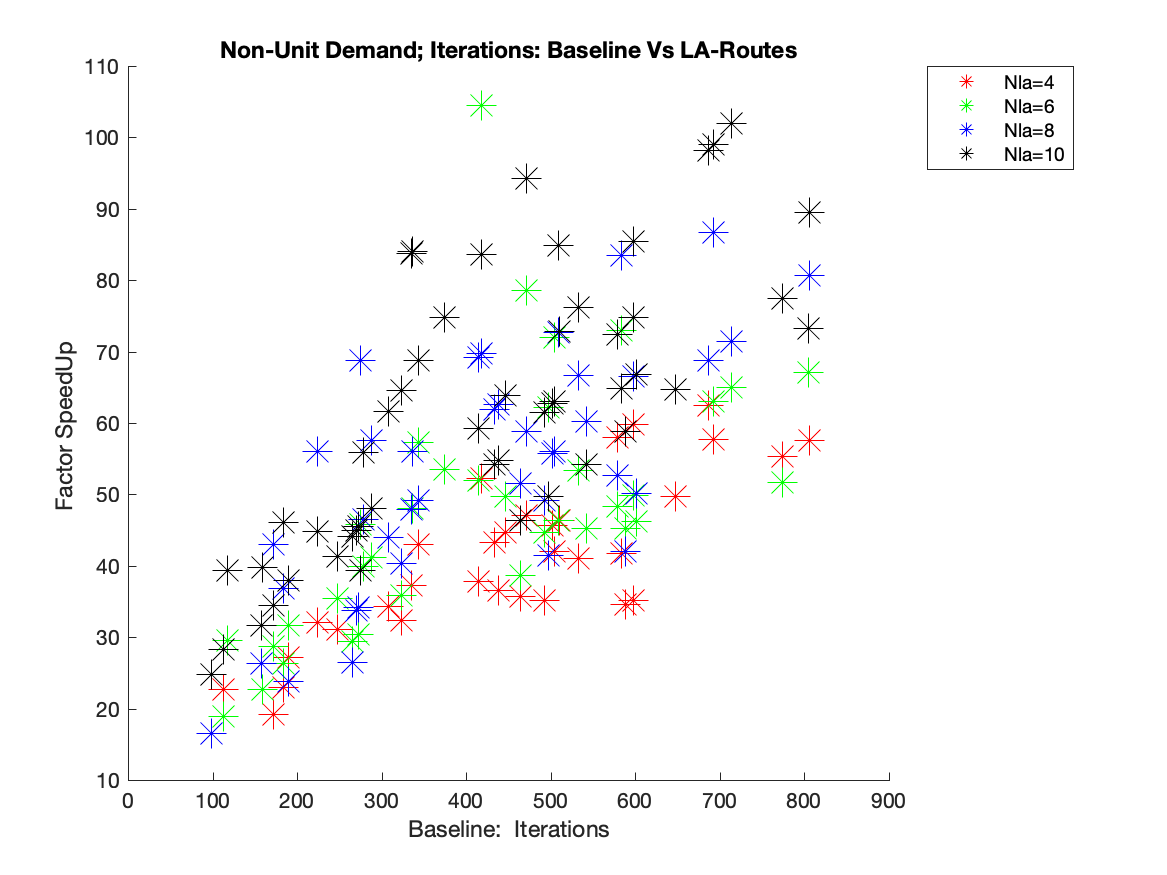}\\
\includegraphics[width=0.45\linewidth]{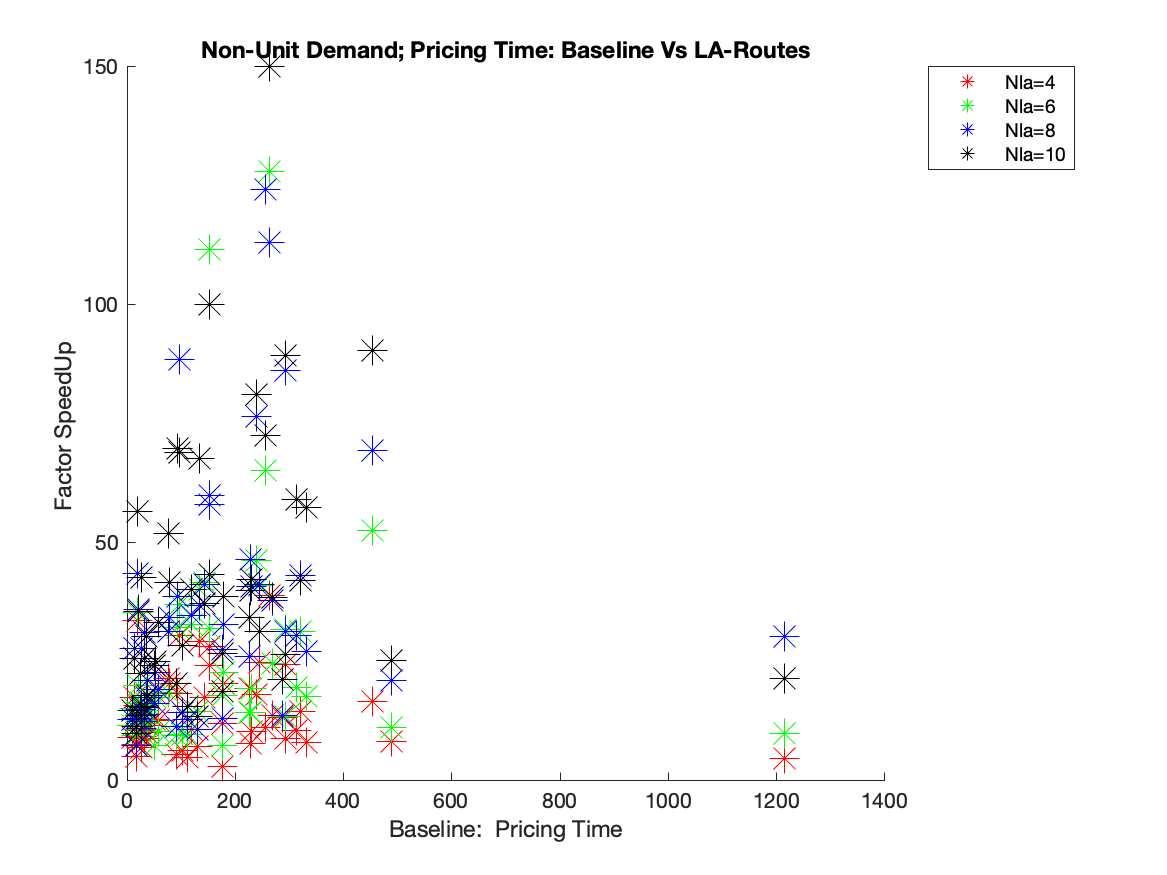}
\includegraphics[width=0.45\linewidth]{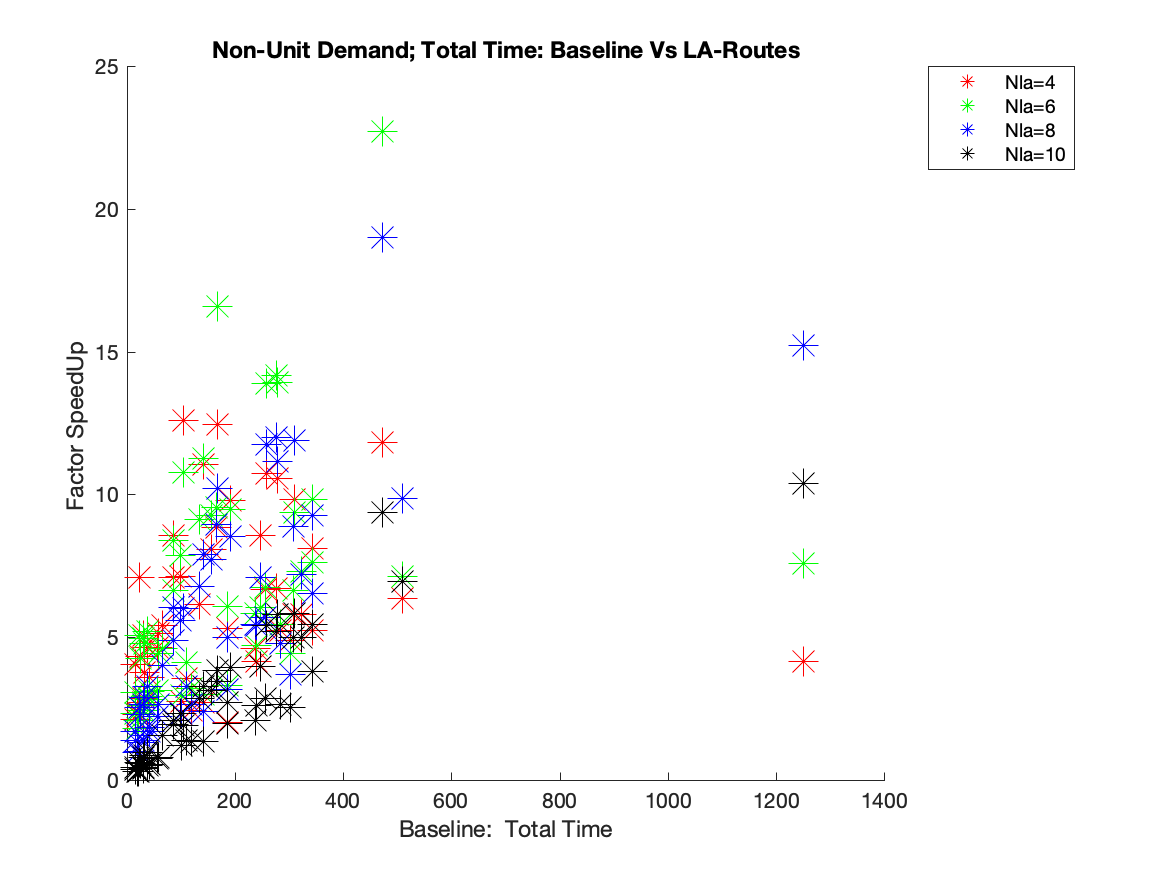}\\

	\caption{In these plots, we the factor speed up relative to the baseline in terms of iterations, pricing time and total time over our data set of non-unit demand problems 
    } 
	\label{fig_unit_3ab}
\end{figure}

We now consider the efficiency of our solver for the master problem without LA-SRI relative to a baseline solver using Decremental State Space Relaxation (DSSR) \citep{righini2009decremental} for pricing and no stabilization.   
Non-unit demand data is generated with ten problem instances for each class where a class is defined by $|N|,d_0$ being in $(20,20),(20,30),(20,40),(30,20),(30,30)$; where each customer has integer demand in range $[1,10]$  uniformly distributed. We display timing results in terms of factor speed up for pricing time, total time and iterations in Fig \ref{fig_unit_3ab} and observe speed ups in each relative to the baseline. 

\section{Combining Rounded Capacity Inequalities with LA-SRI}
\label{RCI_Augerat}

%Just before the conclusion of the paper we hint that exploring LA-SRI combined with rounded capacity inequalities (RCI) (sometimes called rounded capacity constraints) should yield additional benefits over including RCI alone. Here we show the relative improvement possible. We can see consistent and sometimes dramatic gap reductions. 

We now consider the value of Rounded capacity constraints \citep{archetti2011column}. RCI are a powerful class of valid inequality that can be efficiently integrated into pricing for vehicle routing problems. RCI do not alter the structure of pricing and can be used jointly with other valid inequalities. Furthermore they are easily integrated into our fast RMP solution and pricing.  

We now describe RCI. For any set of customers $\hat{N}\subseteq N$ let $b_{\hat{N}}$ be a lower bound on the  minimum number of vehicles servicing at least one customer in $\hat{N}$.  Thus the minimal number of vehicles servicing at least one customer in $\hat{N}$ can be no less than $b_{\hat{N}}$ where $b_{\hat{N}}$ can be naively computed as the total demand of the customers in $\hat{N}$ divided by $d_0$ rounding up.  

Thus inequalities RCI constraints  can be written as follows in the context of the CG-RMP.
\begin{align}
\label{valid1}
    \sum_{l \in \Omega}(\max_{u\in \hat{N}}a_{ul})\theta_l\geq  b_{\hat{N}}\\
     b_{\hat{N}}=\lceil \frac{\sum_{u \in \hat{N}}d_u}{d_0} \rceil
\end{align}
Pricing over inequalities of the form \eqref{valid1} is not easily done because it does not decouple in terms of edges.  The number of vehicles servicing customers in $\hat{N}$ can be upper bounded by the number of edges starting inside $\hat{N}$ and outside $\hat{N}$.  Using $x_{uv}=\sum_{l \in \Omega}\sum_{d\in D}\theta_l a_{uvdl}$ to denote the number of selected edges starting at $u$ and ending at $v$ we write the following valid inequality as follows. 
\begin{align}
\label{valid2}
    \sum_{u \in \hat{N}}\sum_{v \in N^+-\hat{N}}x_{uv}\geq  b_{\hat{N}}
\end{align}
When inequalities are of the form in \eqref{valid2} are added to the standard set cover RMP pricing can be done since the dual variables become associated with edges.  A tighter inequality than \eqref{valid2} (for a given $\hat{N}$) but less tight than \eqref{valid1} is written in terms of LA-arcs below.  
\begin{subequations}
\label{valid_3}
\begin{align}
    \sum_{l \in \Omega_R}\sum_{\substack{y\in Y^l\\ p \in \hat{\Omega}^l_y\\ p=(u,v,N_p)}}a_{\hat{N}p}x^l_{p}\geq b_{\hat{N}}\\
    a_{\hat{N}p}=[0<|\hat{N}\cap (u\cup N_p)|][v \notin \hat{N}]
\end{align}
\end{subequations}
Observe that \eqref{valid_3} is tighter than \eqref{valid2} since entering and existing the set $\hat{N}$ in an LA-arc is treated as occurring once not multiple times (if it occurs multiple times).  Inequalities of the form \eqref{valid_3} can be separated by enumerating over a pre-specified list and added in a cutting plane manner to the RMP.  This cost is associated with LA-arcs and hence does not alter the structure of pricing.  

\begin{figure*}[!hbtp]
\begin{minipage}{1\textwidth}
\includegraphics[width=0.6\linewidth]{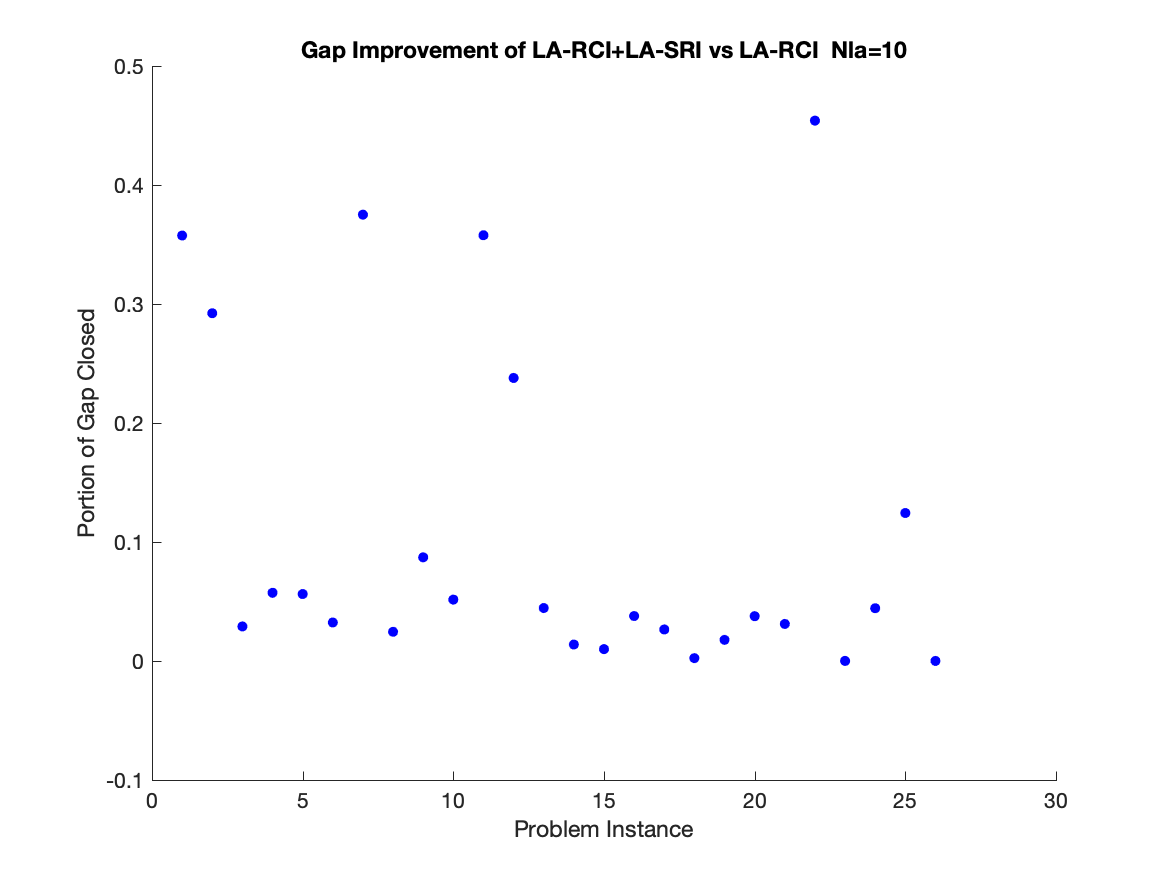}

\caption{The decrease in integrality gap due to LA-SRI over RCI.} 
\label{fig:ROI_Augerat}
 \end{minipage}
\end{figure*}

We use all subsets of nearby customers in addition to denote the set of LA-RCI. In addition we include on LA-SRI $\hat{N}\leftarrow N$ to enforce that at least $\frac{\sum_{u \in N}d_u}{d_0}$ routes are used. We define the RCI under consideration based on nearby neighbors as any subset of the LA neighbors of any customer (and that customer). We write this definition mathematically as follows:  

\begin{align}
\label{sets_use}
    \hat{N} \subseteq N \quad \exists u \in N \quad \hat{N}\subseteq (N_u \cup u)
\end{align}
%In addition we use one inequality where $\hat{N}\leftarrow N$.  
%
We now demonstrate the value of LA-RCI when used jointly with LA-SRI.  We establish that LA-SRI close much of the integrality gap remaining left by using LA-RCI only.  We use eight La neighbors when computing the RCI sets to consider as described in \eqref{sets_use} and ten LA neighbors and LA-SRI selection equal to b ($|N_{\delta}|=3,m_{\delta}=2$ and $|N_{\delta}|=4,m_{\delta}=3$). In Fig \ref{fig:ROI_Augerat} we plot the portion of the integrality gap closed  using LA-RCI vs  LA-RCI+LA-SRI on the benchmark data \citep{augerat1995computational} (with divisor=10)for each problem instance. 
%\newpage

%\section{Other Extensions}
%\label{other_ext}

%We now consider additional extensions of our approach for use in future work.  We seek to use time windows in our formulation.  This may involve different LA neighbor settings depending on the time when the vehicle leaves a customer. We also intend to use the original SRI in our formulation when they are violated by the LP solution but the associated LA-SRI are not.
 % In future work we seek to develop an algorithm that constructs the LA neighbor sets for customers in an efficient manner so as to maximally tighten the LP. In this process, we would aim to maximally tighten the LP bound given that the number of LA neighbors per customer does not exceed a user defined value (thus describing a computation budget). 
%In addition we can adapt other classes of valid inequalities to use with LA-arcs as we did SRI such as capacity inequalities \citep{archetti2011column}.
%We shall also exploit the primal block angular (PBA) structure of the equality constraints in  $\Psi(\Omega_R,\Delta_R)$ (which are separable by $l \in \Omega_R$) for efficient solution as is done for PBA LPs in \citep{castro2007interior}. 

%\end{APPENDIX}

\end{document}